\documentclass{article}[11pt]

\usepackage{algorithm}
\usepackage{algorithmicx}
\usepackage{algpseudocode}

\usepackage{amsmath}
\usepackage{amssymb}
\usepackage{graphicx}
\usepackage{epic,eepic}
\usepackage{comment}
\usepackage{amstext,amsthm,amssymb,amsmath}
\usepackage{epsfig,subfigure,epstopdf}
\usepackage{graphicx}
\usepackage{subfigure}
\usepackage{wrapfig}
\usepackage{fullpage}
\usepackage{color}
\usepackage{multirow}
\usepackage{tabulary}
\usepackage{booktabs}
\usepackage{enumerate}
\usepackage{color}

\renewcommand{\div}{\mathop{\rm div}\nolimits}

\theoremstyle{definition}
\newtheorem{theorem}{Theorem}[section]
\newtheorem{lemma}{Lemma}[section]

\usepackage[utf8]{inputenc}
\usepackage[english]{babel}
\usepackage{hyperref}

\usepackage{nomencl}

\newtheorem{exmp}{Example}[section]

\newcommand{\prnt}[1]{\left( #1 \right)}

\newcommand{\norm}[1]{\left\|#1\right\|}
\newcommand{\abs}[1]{\left |#1\right |}

\newcommand{\normL}[2]{\norm{#1}_{L^2\prnt{#2}}}

\DeclareMathOperator{\argmin}{argmin}
\DeclareMathOperator{\spa}{span}
\DeclareMathOperator{\supp}{supp}

\usepackage{babel}

\begin{document}

\title{Iterative Oversampling Technique for Constraint Energy Minimizing Generalized Multiscale Finite Element Method in the Mixed Formulation}
\author{
Siu Wun Cheung\thanks{Department of Mathematics, Texas A\&M University, College Station, TX 77843, USA 
(\texttt{tonycsw2905@math.tamu.edu})}
\and
Eric Chung\thanks{Department of Mathematics, The Chinese University of Hong Kong, Shatin, New Territories, Hong Kong SAR, China 
(\texttt{tschung@math.cuhk.edu.hk})}
\and
Yalchin Efendiev\thanks{Department of Mathematics \& Institute for Scientific Computation (ISC), Texas A\&M University,
College Station, TX 77843, USA (\texttt{efendiev@math.tamu.edu})}
\and
Wing Tat Leung\thanks{Department of Mathematics, 
University of California Irvine, Irvine, CA 92697, USA (\texttt{wtleung@uci.edu})}
\and
Sai-Mang Pun\thanks{Department of Mathematics, Texas A\&M University, College Station, TX 77843, USA 
(\texttt{smpun@math.tamu.edu})}
}

\maketitle

\begin{abstract}
In this paper, we develop an iterative scheme to construct multiscale basis functions within the framework of the Constraint Energy Minimizing Generalized Multiscale Finite Element Method (CEM-GMsFEM) for the mixed formulation. The iterative procedure starts with the construction of an energy minimizing snapshot space that can be used for approximating the solution of the model problem. A spectral decomposition is then performed on the snapshot space to form global multiscale space. Under this setting, each global multiscale basis function can be split into a non-decaying and a decaying parts. The non-decaying part of a global basis is localized and it is fixed during the iteration. Then, one can approximate the decaying part via a modified Richardson scheme with an appropriately defined preconditioner. Using this set of iterative-based multiscale basis functions, first-order convergence with respect to the coarse mesh size can be shown if sufficiently many times of iterations with regularization parameter being in an appropriate range are performed. Numerical results are presented to illustrate the effectiveness and efficiency of the proposed computational multiscale method. 
\end{abstract}


\section{Introduction}
Many problems arising from engineering involve heterogeneous materials which have strong contrasts in their physical properties.
In general, one may model these so-called multiscale problems using partial differential 
equations (PDEs) with high-contrast valued multiscale coefficients. 
An important example is Darcy's law describing flow in highly heterogeneous porous media. 
These problems are prohibitively costly to solve when traditional fine-scale solvers are directly applied. 
The direct simulation of multiscale PDEs with accurate resolution can be costly as a relatively fine mesh is required to resolve the coefficients, leading to a prohibitively large number of degrees of freedom, a high percentage of which may be extraneous. 
Therefore, some types of model-order reductions are necessary to avoid high computational cost in simulation. 

These computational challenges have been addressed by the development of 
efficient model reduction techniques and many model reduction techniques have been well explored in existing literature. 
For instance, in upscaling methods \cite{chen2003coupled, durlofsky1991numerical, wu2002analysis, gao2015numerical} which are commonly used, one typically derives upscaling media based on the model problem and solves the resulting upscaled problem globally on a coarser grid. In general, this derivation can be done by solving a class of local cell problems in the coarse elements. Besides upscaling approaches mentioned above, multiscale methods  \cite{efendiev2009multiscale,efendiev2000convergence,hou1997multiscale,hughes1995multiscale} have been widely used to approximate the solution of the multiscale problem. 
In multiscale methods, the solution to the problem is approximated using local basis functions, which are solutions to a class of local problems, which are related to the model, on the coarse grid. 
Moreover, because of the necessity of the mass conservation for velocity fields, many approaches have been proposed to guarantee this property, such as multiscale finite volume methods  \cite{cortinovis2014iterative,hajibeygi2011hierarchical,jenny2003multi,lie2017feature,lunati2004multi}, mixed multiscale finite element methods (MsFEM) \cite{aarnes2004use,aarnes2008mixed,chen2003mixed} and its generalization GMsFEM \cite{chan2016adaptive,chen2016least,chung2015mixedg,chung2016mixed,chung2016adaptive,chung2016mixedwave}, mortar multiscale methods \cite{arbogast2007multiscale,chung2016enriched,peszynska2005mortar,peszynska2002mortar} and various post-processing methods \cite{bush2013application,odsaeter2017postprocessing}.


Among these multiscale methods mentioned above, we focus on the framework of GMsFEM in this work. The GMsFEM in mixed formulation has been developed in \cite{chung2015mixedg,efendiev2013mini} recently and it provides a systematic procedure to construct multiple basis functions for either velocity or pressure in each local patch, which makes these methods different to previous methodology in applications. The computation of velocity basis functions involves a construction of snapshot space and a model reduction via local spectral decomposition to identify appropriate modes to form the multiscale space. The convergence analysis in \cite{chung2015mixedg} addresses a spectral convergence with convergence rate proportional to $\Lambda^{-1}$, where $\Lambda$ is the smallest eigenvalue whose modes are excluded in the multiscale space. 
In \cite{chung2018constraintmixed}, a variation of GMsFEM based on a constraint energy minimization (CEM) strategy for mixed formulation has been developed. 
This approach is inspired by the work on localization \cite{maalqvist2014localization,owhadi2017multigrid,owhadi2014polyharmonic} 
and makes use of the ideas of oversampling to compute multiscale basis functions in oversampled subregions with the satisfaction of an appropriate orthogonality condition, where similar ideas have been applied for various numerical discretization and model problems \cite{chung2018constraint,cheung2018mc,cheung2020wave,chung2020wave,li2019parabolic,cheung2019dg}. 
The method proposed in \cite{chung2018constraintmixed} provides a mass conservative velocity field and allows one to identify some non-local information depending on the inputs of the problem. 
One can show that the CEM-GMsFEM provides a better convergence rate (comparing to the mixed GMsFEM) that is proportional to $H\Lambda^{-1}$ with $H$ the size of coarse mesh if the size of oversampling regions is at least of the logarithmic magnitude of the product of the coarse mesh size and the value of contrast. 


However, in the original framework of CEM-GMsFEM, one needs to construct the multiscale basis functions supported in relatively large oversampling regions in order to guarantee a certain level of accuracy and it leads to a moderately large computational cost in the offline stage. 
In particular, when dealing with the case of high-contrast permeability, one needs to set the oversampling parameters to be logarithm of the value of contrast and it results in a loss of sparsity of the stiffness and mass matrices. 

In this work, we propose an iterative computational scheme to construct multiscale basis functions (for velocity) satisfying the property of CEM to overcome the issue mentioned above and enhance the computational efficiency. 
The proposed method relates to the theory of iterative solvers and subspace decomposition methods \cite{kornhuber2016numerical, kornhuber2018analysis,peterseim2018domain} (see also the discussion in \cite[Remark 2.6]{engwer2019efficient} about the iterative implementation of a class of numerical homogenization methods). 
For the mixed formulation, the construction of multiscale basis functions in this work slightly differs from that of the original CEM-GMsFEM and it starts with a set of localized energy minimizing snapshot functions. 
Conceptually, we decompose the global multiscale basis function with CEM property into a decaying and a non-decaying parts. The non-decaying part is formed by the energy minimizing snapshots thus it is localized and it will be fixed during the iterations. 
Then, starting with a zero initial condition, we approximate the decaying part using an iterative scheme of the type of modified Richardson \cite{richardson1910,saad_book} (or any other iterative methods) with an appropriate designed preconditioner. 
The size of support of the approximated decaying part is proportional to the number of the iterations, which can be freely adjusted by the user. Hence, this (iterative) construction for multiscale basis functions is more flexible than the one in the original CEM-GMsFEM. 
The iterative process will maintain the property of mass conservation and no need for any post-processing techniques. 
The proposed method has an advantage that the marginal computational cost from one iteration to next one is comparatively low. 
This iterative construction also shows some potential to compute the offline CEM basis functions in an adaptive manner to further reduce the cost of computation. 
With careful selection of regularization parameter in the iterative scheme, one can show the first-order convergence rate (with respect to $H$) of the velocity if sufficiently many iteration times, depending only on the coarse mesh and the quantity $\Lambda$ in GMsFEM, are performed in the offline stage. 

The paper is organized as follows. In Section \ref{sec:prelim}, we present some preliminaries of the model problem considered in this work. We also briefly review the framework of the original CEM-GMsFEM. Then, we derive the iterative construction of multiscale basis functions for velocity in Section \ref{sec:iterative}. Next, we provide a complete analysis of the proposed iterative construction in Section \ref{sec:analysis}. In particular, we estimate the condition number of the matrix in the iteration in Lemma \ref{lemma:precond-eig}. The main theoretical results of the sufficient condition of linear convergence reads in Theorem \ref{thm:main}. Several numerical tests are provided in Section \ref{sec:numerics} to demonstrate the performance of the numerical methods based on the iterative scheme. Finally, some concluding remarks are drawn in Section \ref{sec:conclusion}. 











\section{Preliminaries}\label{sec:prelim}
\subsection{Model problem}
In this section, we introduce the model problem in this work. 
Consider a class of high-contrast flow problems in the following mixed formulation over a computational domain $D \subset \mathbb{R}^d$ ($d =2, 3$) as follows:
\begin{eqnarray}
\begin{split}
\kappa^{-1} v + \nabla p = 0 & \quad \text{in } D, \\
\nabla \cdot v = f & \quad \text{in } D, \\
v \cdot \mathbf{n} = 0 & \quad \text{on } \partial D, \\
\int_D p ~ dx = 0.
\end{split}
\label{eqn:flow}
\end{eqnarray}
Here, $\mathbf{n}$ is the outward unit normal vector field on the boundary $\partial D$. Note that the source function $f \in L^2(D)$ satisfies the following compatibility condition:
$$ \int_D f ~ dx = 0.$$
In this work, we assume that the function $\kappa : D \to \mathbb{R}$ is a heterogeneous coefficient of high contrast. In particular, there are two constants $\kappa_{\min}$ and $\kappa_{\max}$ such that $0< \kappa_{\min} \leq \kappa(x) \leq \kappa_{\max}$ for almost every $x \in D$ and 
$\kappa_{\max}\kappa_{\min}^{-1} \gg 1$. 
Denote $V := H (\text{div};D)$, $Q:= L^2(D)$, $V_0 := \{ v \in V: v \cdot \mathbf{n} = 0 \text{ on } \partial D\}$, and $Q_0 := \{ q \in Q: \int_D q ~ dx = 0\}$. To numerically solve this problem, we consider the following variational formulation: find $(u , p) \in V_0 \times Q_0$ such that 
\begin{eqnarray}
\begin{split}
a(u,v) - b(v,p) &= 0 & \quad \forall v \in V_0, \\
b(u,q) &= (f,q) & \quad \forall q \in Q, 
\end{split}
\label{eqn:model_var}
\end{eqnarray}
where the bilinear forms $a(\cdot,\cdot)$, $b(\cdot,\cdot)$, and $(\cdot,\cdot)$ are defined as follows:
$$ a(v,w) := \int_D \kappa^{-1} v \cdot w ~ dx, \quad b(w,q) :=  \int_D q ~ \nabla \cdot w~ dx, \quad \text{and} \quad (p,q) := \int_D pq ~ dx$$
for all $v, w \in V$ and $p, q \in Q$. We remark that the following inf-sup condition holds: for all $q \in Q_0$, there is a constant $c>0$, which is independent to $\kappa$, such that 
$$ c \normL{q}{D} \leq \sup_{v \in V_0} \frac{b(v,q)}{\norm{v}_{H(\text{div};D)}}.$$

\noindent
Next, we briefly introduce the notions of coarse and fine grids.
Denote by $\mathcal{T}^H$ a coarse-grid partition of the domain $D$ with mesh size $H$ and by $\mathcal{T}^h$ a fine-grid partition of the domain of $D$ with mesh size $h$ with $0 < h \ll H < 1$. 
We assume that each coarse-grid element $K$ in the coarse-grid partition $\mathcal{T}^H$ contains a connection union of some fine-grid elements. Also, we assume that the fine-grid partition $\mathcal{T}^h$ is sufficiently fine to resolve the small-scale information of the solution.
Let $\{K_i \}_{i = 1}^N$ be the set of coarse-grid elements in the coarse-grid partition $\mathcal{T}^H$ and $N$ be the total number of coarse-grid elements. We denote by $\{ x_j \}_{j=1}^{N_c}$ the set of coarse-grid nodes in the coarse-grid partition $\mathcal{T}^H$ and denote its cardinality by $N_c$. Furthermore, we denote the set of all faces of the coarse-grid partition as $\mathcal{E}^H := \{ E_\ell \}_{\ell=1}^{N_e}$ with cardinality $N_e$.  

\subsection{The GMsFEM with constraint energy minimization}
In this section, we outline the framework of the mixed GMsFEM \cite{chung2016adaptive, chung2015mixedg, efendiev2013generalized} with the setting of constraint energy minimization (CEM) \cite{chung2018constraintmixed}.  
The general procedure of the CEM-GMsFEM can be summarized by the following component modules: (i) perform (local) spectral decomposition; and (ii) find the CEM basis functions. 
Throughout the paper, we define $V_0(\Omega)$ and $Q(\Omega)$ to be the restriction of $V_0$ and $Q$ on the subset $\Omega \subset D$, respectively. 


\subsubsection*{\underline{Perform local spectral decomposition}}
Let $K_i \in \mathcal{T}^H$ be a coarse element. 
We consider the local eigenvalue problem over the coarse element $K_i$ as follows: find $(\phi_j^{(i)}, p_j^{(i)}) \in V_0(K_i) \times Q(K_i)$ and $\lambda_j^{(i)} \in \mathbb{R}$ such that 
\begin{eqnarray}\label{eqn:spectral}
\begin{split}
a(\phi_{j}^{(i)}, v) - b(v, p_j^{(i)}) &=  0 & \quad \forall v \in V_0(K_i), \\
b(\phi_{j}^{(i)}, q) & = \lambda_j^{(i)} s_i ( p_j^{(i)},q)& \quad \forall q \in Q(K_i). 
\end{split}
\end{eqnarray}
where the bilinear form $s_i(\cdot,\cdot)$ is defined as follows:
$$ s_i(p,q) := \int_{K_i} \tilde \kappa pq ~ dx, \quad \tilde \kappa := \kappa \sum_{j=1}^{N_c} \abs{\nabla \chi_j}^2.$$
Here, $\{ \chi_j \}_{j=1}^{N_c}$ is a set of standard multiscale basis functions satisfying the property of partition of unity. 
Specifically, the function $\chi_j$ satisfies the following system
\begin{eqnarray*}
    -\nabla \cdot \left (\kappa(x) \nabla \chi_j \right) = &  0 
                    &\quad \text{in }  K \subset \omega_j, \\
    \chi_j =  & g_j & \quad \text{on } \partial K, \\
    \chi_j = & 0 & \quad \text{on } \partial \omega_j,
\end{eqnarray*}
for all coarse elements $K \subset \omega_j$, where $g_j$ is a linear and continuous function defined on $\partial K$. 
We remark that the definition of $\tilde \kappa$ is motivated by the analysis. 
Next, we arrange the eigenvalues of \eqref{eqn:spectral} in ascending order (i.e. $0 = \lambda_1^{(i)} \leq \lambda_2^{(i)} \leq \cdots$) and we select the first $L_i$ eigenfunctions $\{ p_j^{(i)} \}_{j=1}^{L_i}$ corresponding to the small eigenvalues in order to form an auxiliary space 
$$ Q_{\text{aux}} := \bigoplus_{i=1}^N Q_{\text{aux}}^{(i)}, \quad \text{where }~ Q_{\text{aux}}^{(i)} := \text{span} \left \{ p_j^{(i)}:  j =1, 2, \cdots, L_i  \right \}.$$
Without loss of generality, we assume that $s_i( p_j^{(i)}, p_j^{(i)})= 1$ for all $i \in \{ 1, \cdots, N\}$ and $j \in \{ 1, \cdots, L_i \}$. 
We remark that this auxiliary space $Q_{\text{aux}}$ is used to approximate the pressure $p$ in \eqref{eqn:model_var}. We also define the interpolation operator $\pi : Q \to Q_{\text{aux}}$ as follows:
$$ \pi q = \pi (q) := \sum_{i=1}^N \sum_{j=1}^{L_i} s_i(p_j^{(i)}, q) p_j^{(i)} \quad \forall q\in Q.$$
Note that $\pi$ is the $L^2$-projection of $Q$ to $Q_{\text{aux}}$ with respect to the inner product $s(\cdot, \cdot) := \sum_{i=1}^N s_i(\cdot,\cdot)$. 

\subsubsection*{\underline{Find CEM bases}}
Next, we construct multiscale basis functions satisfying the constraint energy minimization. To be specific, we construct the multiscale basis functions $\psi_{j,{ms}}^{(i)}$ using the auxiliary space $Q_{\text{aux}}^{(i)}$. 
For each $p_j^{(i)} \in Q_{\text{aux}}^{(i)}$, we consider the following system of equations: find $(\psi_j^{(i)}, q_j^{(i)}, \mu_j^{(i)}) \in V_0 \times Q \times Q_{\text{aux}}$ such that 
\begin{eqnarray}
\begin{split}
a(\psi_j^{(i)}, v) - b(v, q_j^{(i)}) = & 0 & \quad \forall v \in V_0, \\
b(\psi_j^{(i)}, q) - s(\mu_j^{(i)},\pi q)  = & 0 & \quad \forall q \in Q, \\
s(\pi q_j^{(i)}, \gamma)  = & s(p_j^{(i)}, \gamma) & \quad \forall \gamma \in Q_{\text{aux}}. 
\end{split}
\label{eqn:cem_1}
\end{eqnarray}
We define $V_{\text{glo}} := \text{span} \{ \psi_j^{(i)}: i = 1, \cdots, N, ~ j = 1, \cdots, L_i \}$. The multiscale space $V_{\text{glo}}$ provides a good approximation (in the sense of Galerkin projection) of the solution $(u,p)$ to the problem \eqref{eqn:model_var}. We define the global solution $(u_{\text{glo}}, p_{\text{glo}}) \in V_{\text{glo}} \times Q_{\text{aux}}$ such that 
\begin{eqnarray}\label{eqn:glo_sol}
\begin{split}
a(u_{\text{glo}}, v) - b(v, p_{\text{glo}}) &=  0 & \quad \forall v \in V_{\text{glo}}, \\
b(u_{\text{glo}}, q) & = (f,q) & \quad \forall q \in Q_{\text{aux}}. 
\end{split}
\end{eqnarray}
However, the function $\psi_j^{(i)}$ usually has a global support, which is computationally infeasible. One can obtain a set of local multiscale basis functions $\left \{ \psi_{j, \text{ms}}^{(i)} \right \}$ by solving a class of the problems \eqref{eqn:cem_1} on each oversampled region of the coarse element, instead of the whole domain. With this set of localized basis functions, the results in \cite{chung2018constraintmixed} show that one can achieve first-order convergence rate (with respect to the coarse mesh size $H$) independent of the contrast provided that sufficiently large oversampling size is considered. 

\section{Iterative construction of multiscale basis functions}\label{sec:iterative}

In this section, we propose an alternative approach to construct the multiscale basis functions satisfying the property of constraint energy minimization. 
Here, the underlying construction is performed based on an iterative process. 
In order to obtain local basis functions, the iterative method is required to keep the support of the basis function in the next iteration is within one or few coarse layers larger than that of the previous iteration. 

\subsection{Construction of offline space}
For each $E_\ell \in \mathcal{E}^H$, we define the coarse neighborhood $\omega_i$ to be $\omega_\ell := \bigcup \{ K \in \mathcal{T}^H: E_\ell \subset K\}$. We denote the set of fine-grid faces (in $\mathcal{T}^h$) lying on the coarse-grid face $E_\ell$ as $\mathcal{E}_h(\omega_\ell)$ with $\abs{\mathcal{E}_h (\omega_\ell) } := J_\ell$. 
Denote $Q_0(\omega_\ell)$ the subspace of $Q(\omega_\ell)$ with zero average piecewise on $K \subset \omega_\ell$.  
A snapshot pair of functions $(\psi_{\text{snap}}^{\ell,j}, p_{\text{snap}}^{\ell,j}) \in V_0(\omega_\ell) \times Q_0(\omega_\ell)$ is the solution of the following system of equations: 
\begin{equation}\label{eqn:def-snap}
\begin{aligned}
\kappa^{-1} \psi_{\text{snap}}^{\ell,j} + \nabla p_{\text{snap}}^{\ell,j}  = &0&\quad \text{in any } K\subset \omega_\ell,\\
\nabla \cdot \psi_{\text{snap}}^{\ell,j} = & \alpha_K I_K & \quad \text{in any } K \subset \omega_\ell,\\
\psi_{\text{snap}}^{\ell,j} \cdot \mathbf{m}_\ell = & \delta_j^\ell  & \quad \text{on }~ E_\ell, \\
\psi_{\text{snap}}^{\ell,j} \cdot \mathbf{n}_\ell = & 0 & \quad \text{on }~ \partial \omega_\ell \setminus E_\ell. 
\end{aligned}
\end{equation}
Here, $I_K$ is the indicator function on $K$, 
$\mathbf{m}_\ell$ is a fixed unit vector field orthogonal to $E_\ell$, and $\mathbf{n}_\ell$ is the unit outward normal vector field of the boundary $\partial \omega_\ell$. 
The function $\delta_j^\ell$ is defined on each fine edge $e_k \in \mathcal{E}_h(\omega_{\ell})$ such that $\delta_j^\ell (e_k) := \abs{e_k}^{-1} \delta_{jk}$ for any $e_k \in \mathcal{E}_h(\omega_\ell)$. 
We remark that $\alpha_K$ is a constant satisfying the following compatibility condition: 
$$ \alpha_K = \vert K \vert^{-1} \int_{E_\ell} \delta_j^\ell ~ dS \quad \forall K \subset \omega_\ell.$$
The snapshot space $V_{\text{snap}} \subset V$ is defined as follows:
$$ V_{\text{snap}} := \bigoplus_{\ell=1}^{N_e} V_{\text{snap}}^{(\ell)}, \quad \text{where} \quad 
V_{\text{snap}}^{(\ell)} := \text{span} \{ \psi_{\text{snap}}^{\ell,j} : j = 1, 2, \cdots, J_\ell \}.$$
The snapshot space $V_{\text{snap}}$ takes care of the flux effects, and its complement $\tilde{V}_0$ in $V_0$, i.e. $V_0 = V_{\text{snap}} \oplus \tilde{V}_0$, is given by 
$$ \tilde{V}_0 = \left\{ v \in V_0: v \cdot\mathbf{m}_\ell = 0 \quad  \forall E_\ell \in \mathcal{E}^H \right\}. $$
We also denote by $Q_H$ the subspace of coarse-scale piecewise constant functions in $Q_0$, i.e. 
$$
Q_H = \left\{ q \in Q_0 : q\vert_K \in P^0(K) \quad \forall K \in \mathcal{T}^H \right\},  
$$
where $P^0(K)$ denotes the space of piecewise constant functions on $K$.

Next, we decompose the local snapshot space into two spaces: $V_{\text{snap},1}^{(\ell)}$ and $V_{\text{snap},2}^{(\ell)}$. The first space $V_{\text{snap},1}^{(\ell)}$ is defined to be the span of the basis functions $\{\varphi_{1}^{\ell,j}\}_{j=1}^{\mathcal{J}_{\ell,1}}$ such that $\varphi_{1}^{\ell,j} \cdot \mathbf{m}_{\ell} =1$ on $E_{\ell}$. The second space $V_{\text{snap},2}^{(\ell)}$ contains all basis functions with average normal flux equal to zero, namely, 
$$ V_{\text{snap},2}^{(\ell)} := \text{span} \left \{ v \in V_{\text{snap}}^{(\ell)} : \int_{E_\ell} v \cdot \mathbf{m}_{\ell} = 0 \right \}.$$
Then, we can define the global spaces $V_{\text{snap},1} := \bigoplus_{\ell = 1}^{N_e} V_{\text{snap},1}^{(\ell)}$ and $V_{\text{snap},2} := \bigoplus_{\ell = 1}^{N_e} V_{\text{snap},2}^{(\ell)}$. 

Now, we discuss the construction of the offline space by performing a spectral decomposition of the local snapshot space to select some dominant modes. Recall that $\omega_{\ell}$ is the union of all coarse blocks which share the coarse edge $E_{\ell}$. We define an operator $\mathcal{H}_{\ell}: V_{\text{snap},2}^{(\ell)} \to \tilde V_{\text{snap},2}^{(\ell)} := \text{span} \{ v : v = \left . w \right |_{\omega_\ell} ~ \text{for } w \in V_{\text{snap}}\}$ such that 
$$ \mathcal{H}_\ell (v) := \argmin \left \{ \norm{w}_a : w \in \tilde V_{\text{snap},2}^{(\ell)} ~ \text{and} ~ w = v ~ \text{on} ~ E_\ell \right \}.$$
Then, we define a spectral problem by finding $(\varphi_{2}^{\ell,j}, \sigma_j^{(\ell)}) \in V_{\text{snap},2}^{(\ell)} \times \mathbb{R}$ satisfying 
\begin{eqnarray}
\int_{\omega_\ell} \kappa^{-1}  \mathcal{H}_\ell (\varphi_{2}^{\ell,j}) \cdot \mathcal{H}_\ell (v) ~ dx = \sigma_j^{(\ell)} \int_{\omega_\ell} \kappa^{-1} \varphi_{2}^{\ell,j} \cdot v~ dx \quad \forall v \in V_{\text{snap},2}^{(\ell)}.
\label{eqn:spectral-h}
\end{eqnarray}
Assume that the eigenvalues obtained from \eqref{eqn:spectral-h} are in ascending order and we define the simplified local offline space as $V_{\text{sms}}^{(\ell)} := \spa \{ \varphi_{2}^{\ell,j} : j = 1,\cdots, \mathcal{J}_{\ell,2} \}$ with $\mathcal{J}_{\ell,2} \in \mathbb{N}^+$. Moreover, we denote $\tilde V_{\text{sms}}^{(\ell)}$ as the orthogonal complement of $V_{\text{sms}}^{(\ell)}$ in $V_{\text{snap},2}^{(\ell)}$. That is, we have the relation $V_{\text{snap}}^{(\ell)} = V_{\text{snap,1}}^{(\ell)} \oplus V_{\text{sms}}^{(\ell)} \oplus \tilde V_{\text{sms}}^{(\ell)}$ with $V_{\text{sms}}^{(\ell)}\perp_{a}\tilde{V}_{\text{sms}}^{(\ell)}$. We denote $V_{\text{sms}} := \bigoplus_{\ell = 1}^{N_e} V_{\text{sms}}^{(\ell)}$ and $\tilde V_{\text{sms}} := \bigoplus_{\ell = 1}^{N_e} \tilde V_{\text{sms}}^{(\ell)}$. Note that $\tilde V_{\text{sms}} \neq V_{\text{sms}}^{\perp}$. 
For each coarse edge $E_\ell \in \mathcal{E}^H$, we enumerate the basis functions 
$\{\varphi_{1}^{\ell,j}\}_{j=1}^{\mathcal{J}_{\ell,1}}$ in $V_{\text{snap},1}^{(\ell)}$ and 
$\{\varphi_{2}^{\ell,j}\}_{j=1}^{\mathcal{J}_{\ell,2}}$ in $V_{\text{sms}}^{(\ell)}$ 
as $\{\varphi_j^{(\ell)}\}_{j=1}^{\mathcal{J}_{\ell}}$,  
which constitutes a basis for $V_{\text{snap},1}^{(\ell)} \oplus V_{\text{sms}}^{(\ell)}$ 
containing $\mathcal{J}_{\ell} = \mathcal{J}_{\ell,1} + \mathcal{J}_{\ell,2}$ basis functions. 
For each basis function $\varphi_j^{(\ell)} \in V^{(\ell)}_{\text{snap},1} \oplus V^{(\ell)}_{\text{sms}}$, 
we define a corrector function $\tilde \psi_{\text{ms},j}^{(\ell)} \in \tilde V_{\text{sms}}$ such that it solves the following equation: 
\begin{eqnarray}\label{eqn:main-iterative}
a\left (\tilde \psi_{\text{ms},j}^{(\ell)}, v \right) = - a\left (\varphi_j^{(\ell)}, v\right) \quad \forall v \in \tilde V_{\text{sms}}.
\end{eqnarray}
Then the element $\psi^{(\ell)}_{\text{ms},j} \in V_{\text{snap}}$, defined by 
\begin{equation}\label{eqn:glo_basis}
\psi^{(\ell)}_{\text{ms},j}:=\tilde \psi^{(\ell)}_{\text{ms},j} + \varphi^{(\ell)}_j,
\end{equation}
refers to a global multiscale basis function. 
The global multiscale space is then defined by 
$$V_{\text{glo}} :=\text{span}\left \{\psi^{(\ell)}_{\text{ms},j} \right \} \subset V_{\text{snap}}.$$
It will be shown in our analysis that this set of basis functions $\left \{ \psi_{\text{ms},j}^{(\ell)} \right \}$ provides similar approximability as the CEM basis functions defined in \eqref{eqn:cem_1}. In this work, we use an iterative method for constructing the basis functions $\tilde \psi_{\text{ms}, j}^{(\ell), k}$ that approximates the basis function $\tilde \psi_{\text{ms},j}^{(\ell)}$ in $k$ steps of iteration. 
The global multiscale model for \eqref{eqn:model_var} is then given by: find 
$(u_{\text{glo}}, p_{\text{glo}}) \in V_{\text{glo}} \times Q_H$ such that 
\begin{eqnarray}\label{eqn:glo_msm}
\begin{split}
a(u_{\text{glo}},v)-b(v,p_{\text{glo}}) & =0& \quad \forall v\in V_{\text{glo}}\\
b(u_{\text{glo}},q) & =(f,q)& \quad \forall q\in Q_{H}. 
\end{split}
\end{eqnarray}

\subsection{Derivation of iterative oversampling scheme}

\noindent
In this section, we discuss the construction of iterative scheme for the basis functions.  By such iterative scheme, we construct $\tilde \psi_{\text{ms}, j}^{(\ell), k}$ that approximates the basis function $\tilde \psi_j^{(\ell)}$ defined in \eqref{eqn:main-iterative} in $k$ steps of iteration. We require that the iterative scheme satisfies two general assumptions: (i) After one level of iteration, the support of the basis function will only be slightly enlarged by one coarse layer; and (ii) the iterative basis function possesses a property of exponential decay outside the support of the basis functions. Following the analysis of mixed CEM-GMsFEM in \cite{chung2018constraintmixed}, this leads to convergence to the exact solution to the problem \eqref{eqn:model_var}. 

\noindent 
We aim to construct a sequence of functions $\left \{ \tilde \psi_{\text{ms}, j}^{(\ell), k} \right \}_{k=1}^\infty$ to approximate the multiscale basis function $\tilde \psi_j^{(\ell)}$. 
For a given domain $\omega \subset D$ and $k \in \mathbb{N}$, we define $\omega^{+,k}$ to be the oversampling region such that 
$$ \omega^{+,k} := \left \{ 
\begin{array}{ll} 
\omega & \text{if } k = 0, \\
\bigcup \{ K : \omega^{+, k-1} \cap K \neq \emptyset \} & \text{if } k \geq 1.
\end{array} \right .$$
First, we assume that the initial guess $\tilde \psi_{\text{ms}, j}^{(\ell),0}$ is zero. In the $k$-th level of iteration, assuming we have already constructed the basis function $\tilde \psi_{\text{ms}, j}^{(\ell),k-1}$ from the previous level, we find $\tilde \psi_{\text{ms}, j}^{(\ell),k} \in \tilde V_{\text{sms}}$ such that 
\begin{eqnarray}\label{eqn:iterative-sch-1}
\tilde\psi_{\text{ms}, j}^{(\ell),k} = \tilde\psi_{\text{ms}, j}^{(\ell),k-1} + \tau \sum_{s\in \mathcal{I}_k} \eta_{\ell, j}^{k,s}
\end{eqnarray}
where $\mathcal{I}_k := \{ s \in \mathbb{N}: \omega_s \subset \omega_{\ell}^{+,k}\}$, $\tau$ is a regularization parameter, and $\eta_{\ell, j}^{k,s} \in \tilde V_{\text{sms}}^{(s)}$ solves the following equation: 
\begin{eqnarray}\label{eqn:iterative-sch-2}
a\left ( \eta_{\ell, j}^{k,s} , v \right ) = -a\left (\tilde\psi_{\text{ms}, j}^{(\ell), k-1} + \varphi_j^{(\ell)}, v\right) \quad \forall v \in \tilde V_{\text{sms}}^{(s)}
\end{eqnarray}
for any $ s\in \mathcal{I}_k$. Clearly, we have $\supp\left (\sum_{s \in \mathcal{I}_k} \eta_{\ell,j}^{k,s}\right) \subset \left ( \supp(\tilde\psi_{\text{ms},j}^{(\ell),k-1})\right )^{+,1}$ and  $\supp(\tilde\psi_{\text{ms},j}^{(\ell),k}) \subset \omega_{\ell}^{+,k}$ inductively.
To simplify notation, we use the single-index notation to represent the multiscale space $\tilde V_{\text{sms}}$. In particular, we write 
$$\tilde V_{\text{sms}} = \spa\{ \theta_1 , \cdots, \theta_{\mathcal{M}} \} \quad \text{with} \quad \mathcal{M} := \text{dim}(\tilde V_{\text{sms}}).$$

The remaining of this section will be devoted to transforming the equation \eqref{eqn:main-iterative} into a matrix form and discussing the relations of our iterative scheme \eqref{eqn:iterative-sch-1} to a well-established iterative method for solving linear systems, namely the modified Richardson iteration \cite{richardson1910,saad_book}. That is, we find the vector $\tilde \Psi_j^{(\ell)} \in \mathbb{R}^\mathcal{M}$ representing the coefficients in the expansion of $\tilde \psi_j^{(\ell)}$ using the basis functions in $\tilde V_{\text{sms}}$ such that 
\begin{eqnarray}\label{eqn:main-iterative-mat}
A \tilde \Psi_j^{(\ell)} = y_j^{(\ell)}, 
\end{eqnarray}
where $A := \left ( a(\theta_j, \theta_i) \right ) \in \mathbb{R}^{\mathcal{M} \times \mathcal{M}}$ and $y_j^{(\ell)} := -\left ( a(\varphi_j^{(\ell)},\theta_1) ~ \cdots ~ a(\varphi_j^{(\ell)},\theta_\mathcal{M}) \right )^T \in \mathbb{R}^{\mathcal{M}}$. 
We denote $A_{\ell_1,\ell_2}$ the (block) matrix representation of the bilinear form $a(\cdot,\cdot)$ 
on $\tilde V^{(\ell_1)}_{\text{sms}} \times \tilde V^{(\ell_2)}_{\text{sms}}$ for arbitrary coarse edges $E_{\ell_1}, E_{\ell_2} \in \mathcal{E}^H$. 
One can obtain the global matrix $A$ of the bilinear form $a(\cdot,\cdot)$ on $\tilde V_{\text{sms}} \times \tilde V_{\text{sms}}$ 
by assembling the submatrices, $A_{\ell_i, \ell_j}$, i.e. 
$$
A=\left[
\begin{array}{cccc}
A_{1,1} & A_{1,2} & \cdots & A_{1,N_e}\\
A_{2,1} & A_{2,2} & \cdots & A_{2,N_e}\\
\vdots & \vdots & \ddots & \vdots\\
A_{N_e,1} & A_{N_e,2} & \cdots & A_{N_e,N_e}\\
\end{array}\right]. 
$$
The graph connectivity of local submatrices of $A$ are subject to overlaps of the coarse neighborhoods, i.e. 
$$
\omega_{\ell_1} \cap \omega_{\ell_2} = \emptyset \implies A_{\ell_1,\ell_2} = \mathbf{0}.
$$
In other words, the matrix representation $A$ is block-sparse. 
Moreover, since the bilinear form $a(\cdot,\cdot)$ is symmetric and elliptic, 
the matrix $A$ is symmetric and positive definite. 

We are now going to present the transformation of the iterative scheme \eqref{eqn:iterative-sch-1} into a matrix form. We denote $Z_{\ell, j}^{k} \in \mathbb{R}^{\mathcal{M}}$ is the vector representing the coefficients in the expansion of $\sum_{s\in \mathcal{I}_k}  \eta_{\ell,j}^{k,s}$ using the basis functions in $\tilde V_{\text{sms}}$. From \eqref{eqn:iterative-sch-2}, thanks to the locality of the function $\eta_{\ell,j}^{k,s}$, we have 
$$ A_{\text{pre}} Z_{\ell,j}^{k} = -(A \tilde\Psi_{\text{ms},j}^{(\ell), k-1} - y_j^{(\ell)}) \implies Z_{\ell_j}^{k} = - A_{\text{pre}}^{-1} (A\tilde\Psi_{\text{ms}, j}^{(\ell), k-1} - y_j^{(\ell)}),$$
where $A_{\text{pre}} \in \mathbb{R}^{\mathcal{M} \times \mathcal{M}}$ is the block diagonal part in the assemble of the global matrix $A$, i.e.  
$$ A_{\text{pre}} :=  \left [ \begin{array}{cccc}
A_{1,1} & 0 & \cdots & 0 \\
0 & A_{2,2} & \cdots & 0 \\
\vdots & \vdots & \ddots & \vdots \\
0 & 0 & \cdots & A_{N_e, N_e} 
\end{array}\right ].$$
Therefore, the iterative scheme \eqref{eqn:iterative-sch-1} can be written in the follow matrix form: 
\begin{eqnarray}\label{eqn:iterative-sch-3}
\tilde\Psi_{\text{ms}, j}^{(\ell),k} = \tilde\Psi_{\text{ms}, j}^{(\ell),k-1} + \tau Z_{\ell,j}^{k} = (\mathbf{I}_{\mathcal{M}} - \tau A_{\text{pre}}^{-1} A)\tilde\Psi_{\text{ms}, j}^{(\ell), k-1} -\tau A_{\text{pre}}^{-1}\Phi_j^{(\ell)}, 
\end{eqnarray}
where $\mathbf{I}_{\mathcal{M}} \in \mathbb{R}^{\mathcal{M} \times \mathcal{M}}$ is the identity matrix on $\mathbb{R}^{\mathcal{M}}$. 
We remark that \eqref{eqn:iterative-sch-3} is equivalent to performing a modified Richardson iteration to the preconditioned system of \eqref{eqn:main-iterative-mat} with the preconditioner $A_{\text{pre}}$. 
We denote the element in $\tilde V_{\text{sms}}$ with coefficients $\tilde{\Psi}^{(\ell)}_{\text{ms},j}$ by $\tilde \psi^{(\ell)}_{\text{ms},j}$. 
Then the element $\psi^{(\ell),k}_{\text{ms},j} \in V_{\text{snap}}$, defined by 
\begin{equation}\label{eqn:loc_basis}
\psi^{(\ell),k}_{\text{ms},j}:=\tilde \psi^{(\ell),k}_{\text{ms},j} + \varphi^{(\ell)}_j,
\end{equation}
refers to a localized multiscale basis function approximating the global multiscale basis function $\psi^{(\ell)}_{\text{ms},j}$ in \eqref{eqn:main-iterative}. The localized multiscale space is then defined by 
$$V_{\text{ms}} :=\text{span}\left \{\psi^{(\ell),k}_{\text{ms},j} \right \} \subset V_{\text{snap}}.$$
The localized multiscale model for \eqref{eqn:model_var} is then given by: find 
$(u_{\text{ms}}, p_{\text{ms}}) \in V_{\text{ms}} \times Q_H$ such that 
\begin{eqnarray}\label{eqn:loc_msm}
\begin{split}
a(u_{\text{ms}},v)-b(v,p_{\text{ms}}) & =0 & \quad \forall v\in V_{\text{ms}},\\
b(u_{\text{ms}},q) & =(f,q) & \quad \forall q\in Q_{H}. 
\end{split}
\end{eqnarray}

\section{Analysis}\label{sec:analysis}
In this section, we present some theoretical results of the proposed iterative construction for multiscale basis function satisfying the property of constraint energy minimization. 
We start with introducing some notations which will 
facilitate our discussion. 
To begin with, we define the following $a$-induced weighted $L^2$ norm $\| \cdot \|_a$ on the space $V$: 
\[
\| v \|_a := a(v,v)^{\frac{1}{2}} = \left(\int_D \kappa^{-1} \vert v \vert^2 ~ dx \right)^\frac{1}{2}. 
\]
Throughout this section, we write $a \lesssim b$ is there exists a generic constant $C>0$ such that $a \leq Cb$. 
For any symmetric and positive definite matrix $\mathbf{K}$, 
the norm  $\| \cdot \|_\mathbf{K}$ is defined as 
$\| \mathbf{\Phi} \|_\mathbf{K} = \left(\mathbf{\Phi}^\mathsf{T} \mathbf{K} \mathbf{\Phi}\right)^{1/2}$.

We recall that we have the decomposition $V_0 = V_\text{snap} \oplus \tilde{V}_0$, 
which implies that for any $v \in V_0$, there exists a unique decomposition 
\[
v = v_\text{snap}+ \tilde{v} \text{ where }
v_{\text{snap}} \in V_{\text{snap}}, \,
\tilde{v} \in \tilde{V}_0.
\]
We denote the projection from $v \in V_0$ 
to $\tilde{v} \in \tilde{V}_0$ by $\tilde{\pi}_0$. 
On the other hand, by the construction of the snapshot functions in \eqref{eqn:def-snap}, we have 
$\nabla \cdot V_{\text{snap}} \subset Q_H$. 
The snapshot solution for \eqref{eqn:model_var}, 
denoted by $(u_{\text{snap}}, p_{\text{snap}}) \in V_{\text{snap}} \times Q_H$, 
is defined by
\begin{equation}
\begin{aligned}
a(u_{\text{snap}},v)-b(v,p_{\text{snap}}) & =0 & \forall v\in V_{\text{snap}},\\
b(u_{\text{snap}},q) & =(f,q) & \forall q\in Q_{H}.
\end{aligned}
\label{eqn:snapshot_solution}
\end{equation}
In addition, recall that we have the decomposition of the snapshot space 
\[
V_{\text{snap}} = V_{\text{snap,1}} \oplus V_{\text{snap,2}} = V_{\text{snap,1}} \oplus V_{\text{sms}} \oplus \tilde V_{\text{sms}},
\] 
which implies that for any $v_\text{snap} \in V_{\text{snap}}$, there exists a unique decomposition 
\[
v_{\text{snap}} = v_{\text{snap},1}+v_\text{sms}+\tilde{v}_\text{sms} \text{ where }
v_{\text{snap},1} \in V_{\text{snap},1}, \,
v_\text{sms} \in V_\text{sms}, \,
\tilde{v}_\text{sms} \in \tilde{V}_\text{sms}.
\]
We denote the projection from $v_{\text{snap}} \in V_{\text{snap}}$ 
to $\tilde{v}_\text{sms} \in \tilde{V}_\text{sms}$ by $\tilde{\pi}_{\text{sms}}$. 
Finally, for the pressure space $Q$, we denote by $\pi_{Q_H}$ the coarse-scale 
piecewise $L^2$ projection onto $Q_H$, i.e. $\pi_{Q_H}: Q \to Q_H$ is defined by 
\[
\pi_{Q_H} q = \sum_{K \in \mathcal{T}^H} \vert K \vert^{-1} (q, 1)_K I_K,\quad \forall q \in Q, 
\]
where $(\cdot, \cdot)_K$ denotes the standard $L^2(K)$ inner product, 

The first lemma states a orthogonality property about the decomposition of the snapshot space. 
\begin{lemma}\label{lemma:ortho-1}
For any $(v,q) \in (\tilde{V}_0 \oplus V_{\text{snap},2}) \times Q_H$, we have $b(v,q) = 0$.
\begin{proof}
Let $v \in \tilde{V}_0 \oplus V_{\text{snap},2}$. For any $E_\ell \in \mathcal{E}^H$, we have 
$\int_{E_{\ell}}v \cdot \textbf{m}_{{\ell}}=0$. By divergence theorem, $\int_K \nabla\cdot v=0$ for all $K \in \mathcal{T}^H$. The result follows directly. 
\end{proof}
\end{lemma}

Now we are going to show that the velocity component of the snapshot solution and the global multiscale solution are in fact identical. 
\begin{lemma}\label{lemma:snap-glo}
Let $(u_{\text{snap}},p_{\text{snap}})$ be the solution in \eqref{eqn:snapshot_solution} 
and $(u_{\text{glo}},p_{\text{glo}})$ be the solution in \eqref{eqn:glo_msm}. 
Then we have $u_{\text{snap}} = u_{\text{glo}}$. 
\begin{proof}
For any $v\in\tilde{V}_{\text{sms}}\subset V_{\text{snap},2}$, by Lemma~\ref{lemma:ortho-1}, we infer from the first equation in \eqref{eqn:snapshot_solution} that 
\begin{equation}
a(u_{\text{snap}}, v) = b(v, p_{\text{snap}}) = 0.
\end{equation}
Since $(I - \tilde{\pi}_{\text{sms}})u_{\text{snap}} \in V_{\text{snap},1} \oplus V_{\text{sms}}$, we can uniquely write it into a linear combination
\[
(I - \tilde{\pi}_{\text{sms}})u_{\text{snap}} = 
\sum_{E_\ell \in \mathcal{E}^H} \sum_{j=1}^{\mathcal{J}_\ell} 
c_{j}^{(\ell)} \varphi_{j}^{(\ell)}.
\]
Therefore, using the definition of the global multiscale basis function in \eqref{eqn:main-iterative}, we have 
\[
\begin{split}
a(\tilde{\pi}_{\text{sms}} u_{\text{snap}}, v)
& = -a((I - \tilde{\pi}_{\text{sms}})u_{\text{snap}}, v) \\
& = -\sum_{E_\ell \in \mathcal{E}^H} \sum_{j=1}^{\mathcal{J}_\ell} 
c_{j}^{(\ell)} a(\varphi_{j}^{(\ell)}, v) \\
& = \sum_{E_\ell \in \mathcal{E}^H} \sum_{j=1}^{\mathcal{J}_\ell} 
c_{j}^{(\ell)} a\left (\tilde \psi_{\text{ms},j}^{(\ell)}, v \right) \\
& = a\left (\sum_{E_\ell \in \mathcal{E}^H} \sum_{j=1}^{\mathcal{J}_\ell} 
c_{j}^{(\ell)} \tilde \psi_{\text{ms},j}^{(\ell)}, v \right). 
\end{split}
\]
By taking $v = \tilde{\pi}_{\text{sms}} u_{\text{snap}} - \sum_{E_\ell \in \mathcal{E}^H} \sum_{j=1}^{\mathcal{J}_\ell} c_{j}^{(\ell)} \tilde \psi_{\text{ms},j}^{(\ell)} \in \tilde{V}_{\text{sms}}$, 
we have  
\begin{equation}
\tilde{\pi}_{\text{sms}} u_{\text{snap}} = \sum_{E_\ell \in \mathcal{E}^H} \sum_{j=1}^{\mathcal{J}_\ell} 
c_{j}^{(\ell)} \tilde \psi_{\text{ms},j}^{(\ell)}. 
\end{equation}
Recalling the definition of the global basis functions in \eqref{eqn:glo_basis}, we have 
\begin{equation}
u_{\text{snap}} = \sum_{E_\ell \in \mathcal{E}^H} \sum_{j=1}^{\mathcal{J}_\ell} 
c_{j}^{(\ell)} \psi_{\text{ms},j}^{(\ell)} \in V_{\text{glo}}. 
\end{equation}
Since $V_{\text{glo}}\subset V_{\text{snap}}$, 
subtracting \eqref{eqn:glo_msm} from \eqref{eqn:snapshot_solution}, 
we have
\begin{equation}
\begin{aligned}
a(u_{\text{snap}}-u_{\text{glo}},v)-b(v,p_{\text{snap}}-p_{\text{glo}}) & =0 & \forall v\in V_{\text{glo}},\\
b(u_{\text{snap}}-u_{\text{glo}},q) & =0 & \forall  q\in Q_{H},
\end{aligned}
\end{equation}
and therefore, by putting $v=u_{\text{snap}}-u_{\text{glo}} \in V_{\text{glo}}$ and $q=p_{\text{snap}}-p_{\text{glo}} \in Q_H$, we have $u_{\text{snap}}=u_{\text{glo}}$.
\end{proof}
\end{lemma}

Next, we are going to analyze the error between the weak solution 
and the global multiscale solution. 
For any $K \in \mathcal{T}^H$, we let 
$\tilde{V}_0(K)$ be the restriction of $\tilde{V}_0$ on $K$ and 
$Q_0(K) = \left\{ q \in L^2(K): \; (q,1)_K = 0 \right\}$ 
be the subspace of $L^2(K)$ functions with average zero.
Then we have the following local orthogonality properties.
\begin{lemma}\label{lemma:ortho-2}
Let $K \in \mathcal{T}^H$ and let $q \in Q_0(K)$. 
For any $q' \in Q_H$ and $v \in V_{\text{snap}}$, we have $(q', q)_K = 0$ and $b(v,q) = 0$. 
\begin{proof}
The first result is trivially true from the definition of $Q_H$ and $Q_0(K)$. 
The second result follows directly from the fact that $\nabla \cdot V_{\text{snap}} \subset Q_H$ by the construction of velocity snapshot functions in \eqref{eqn:def-snap}. 
\end{proof}
\end{lemma}

The following lemma states a local stability result for $\tilde{V}$, 
which will be used to derive an estimate 
for the error between the weak solution 
and the global multiscale solution. 
\begin{lemma}\label{lemma:loc_stab}
Let $K \in \mathcal{T}^H$ and $f_K \in L^2(K)$.
Suppose $(u_K,p_K) \in \tilde{V}_{0}(K) \times Q_0(K)$ satisfies 
\begin{equation}
\begin{aligned}
a(u_K,v)-b(v,p_K) & = 0 & \forall v\in \tilde{V}_{0}(K),\\
b(u_K,q) & = (f_K, q)_K & \forall q\in Q_{0}(K).
\end{aligned}
\end{equation}
Then we have 
\[
\|u_K \|_{a(K)}\leq C\min_{x\in K}\{\kappa(x)\}^{-1/2}H\|f_K \|_{L^{2}(K)}.
\]
\begin{proof}
Let $\hat{K}=[0,1]^{d}$ be a reference element and $K=\prod_{i=1}^d [x_{i},x_{i}+H_i] \in \mathcal{T}^H$. 
We define a bijective affine mapping $T_K : \hat{K}\rightarrow K$ by 
\[ T_K(\hat{x}) = 
(x_1, x_2, \ldots, x_d) + J_k \hat{x} 
\quad \forall \hat{x} = (\hat{x}_1, \hat{x}_2, \ldots, \hat{x}_d) \in \hat{K}, \]
where $J_K = \text{diag}(H_1, H_2, \ldots, H_d)$ is the Jacobian matrix of the affine transformation $T_K$, and define the spaces 
\[
\begin{split}
\tilde{V}_0(\hat{K}) & = \{\hat{v}\in H(\text{div}; \hat{K}): \hat{v}= \det(J_K) J_K^{-1} (v\circ T_K) \text{ for some }v\in \tilde{V}_0(K)\}, \\ 
Q_0(\hat{K}) & = \{\hat{q}\in L^{2}(\hat{K}):\;\hat{q}=q\circ T_K \;\text{ for some }q\in Q_0(K)\}.
\end{split}
\]
We denote $u_{\hat{K}}= \det(J_K) J_K^{-1} (u_{K}\circ T_K) \in \tilde{V}_0(\hat{K})$, 
and $p_{\hat{K}} = p_K \circ T_K \in Q_0(\hat{K})$. Then
$(u_{\hat{K}},p_{\hat{K}}) \in \tilde{V}_{0}(\hat{K}) \times Q_0(\hat{K})$ is the solution satisfying
\begin{equation}
\begin{aligned}
a_{\hat{K}}(u_{\hat{K}},v)-b_{\hat{K}}(v,p_{\hat{K}}) & = 0 & \forall v\in \tilde{V}_{0}(\hat{K}),\\
b_{\hat{K}}(u_{\hat{K}},q) & = \det(J_K) (f_K \circ T_K, q)_{\hat{K}} & \forall q\in Q_0(\hat{K}),
\end{aligned}\label{eqn:local_u_0_re}
\end{equation}
where 
$$ a_{\hat{K}}(v,w) := H^{-2} \int_{\hat{K}} (\kappa^{-1} \circ T_K) (J_K v) \cdot (J_K w) ~ dx, \quad 
\text{and} \quad 
b_{\hat{K}}(w,q) :=  \int_{\hat{K}} q ~ \nabla \cdot w~ dx.$$
for all $v, w \in \tilde{V}_0({\hat{K}})$ and $q \in Q_0({\hat{K}})$. 
By the stability of the problem (\ref{eqn:local_u_0_re}), we have
\[
H^{-1} \| (\kappa^{-1/2} \circ T_K) (J_K u_{\hat{K}}) \|_{L^2(\hat{K})}\leq C \min_{x\in K}\{\kappa(x)\}^{-1/2} \left\vert \det(J_K) \right\vert \|f_K\circ T_K\|_{L^{2}(\hat{K})}, 
\]
This implies 
\begin{align*}
\|u_{K}\|_{a(K)}^{2} 
& = \vert \det(J_K) \vert \left\| (\kappa^{-1/2} \circ T_K) \det(J_K)^{-1} (J_K u_{\hat{K}}) \right\|_{L^2(\hat{K})}^2\\ 
& = \vert \det(J_K) \vert^{-1} \left\| (\kappa^{-1/2} \circ T_K) (J_K u_{\hat{K}}) \right\|_{L^2(\hat{K})}^2\\ 
& \leq C \min_{x\in K}\{\kappa(x)\}^{-1} H^2 \left\vert \det(J_K) \right\vert \|f_K\circ T_K\|_{L^{2}(\hat{K})}^2 \\ 
 & \leq C\min_{x\in K}\{\kappa(x)\}^{-1}H^{2}\|f_K\|_{L^{2}(K)}^{2},
\end{align*}
which completes the proof. 
\end{proof}
\end{lemma}

The following lemma states that the error between 
the weak solution and the global multiscale solution 
converges linearly with the coarse mesh size $H$ 
in the $a$-induced norm. 
\begin{lemma}\label{lemma:weak-glo}
Let $(u,p)$ be the solution in \eqref{eqn:model_var} and 
$(u_{\text{glo}},p_{\text{glo}})$ be the solution in \eqref{eqn:glo_msm}. We have
\[
\|u-u_{\text{glo}}\|_{a}\leq C\min_{x\in\Omega}\{\kappa(x)\}^{-1/2}H\|(I - \pi_{Q_H}) f\|_{L^{2}}
\]
\begin{proof}
First, we have $u_{\text{glo}} = u_{\text{snap}}$ from Lemma~\ref{lemma:snap-glo}. 
Subtracting \eqref{eqn:snapshot_solution} from \eqref{eqn:model_var}, we deduce 
\begin{equation}
\label{eqn:weak-snap}
\begin{aligned}
a(u-u_{\text{snap}},v) -b(v,p-p_{\text{snap}}) & = 0 & \forall v\in V_{\text{snap}}, \\
b(u-u_{\text{snap}},q) & =0 & \forall q\in Q_{H}.
\end{aligned}
\end{equation}
By Lemma~\ref{lemma:ortho-1}, for any $q \in Q_H$, since $\tilde{\pi}_0(u) \in \tilde{V}_0$, we have $b(\tilde{\pi}_0(u), q) = 0$, which leaves the second equation of \eqref{eqn:weak-snap} as 
\[
b((I-\tilde{\pi}_0)u-u_{\text{snap}}, q) = 0 \quad \forall q\in Q_{H}.
\]
We note that $(I-\tilde{\pi}_0)u-u_{\text{snap}} \in V_{\text{snap}}$. 
Using that fact that $\nabla \cdot V_{\text{snap}} \subset Q_H$, we take 
$q = \nabla \cdot ((I-\tilde{\pi}_0)u-u_{\text{snap}}) \in Q_H$ to deduce that 
$\nabla \cdot \left((I-\tilde{\pi}_0)u-u_{\text{snap}}\right)= 0$ almost everywhere in $D$, 
and therefore 
\[
b((I-\tilde{\pi}_0)u-u_{\text{snap}}, p - p_{\text{snap}}) 
= \int_{D} (p - p_{\text{snap}}) \nabla \cdot \left((I-\tilde{\pi}_0)u-u_{\text{snap}}\right) = 0.
\]
Taking $v = (I-\tilde{\pi}_0)u-u_{\text{snap}} \in V_{\text{snap}}$ in the 
first equation of \eqref{eqn:weak-snap}, one has 
\[
a\left((I-\tilde{\pi}_0)u-u_{\text{snap}}, (I-\tilde{\pi}_0)u-u_{\text{snap}}\right) 
= -a\left(\tilde{\pi}_0 u,(I-\tilde{\pi}_0)u-u_{\text{snap}}\right), 
\]
which implies 
\begin{equation}\label{eqn:main-stab}
\| u-u_{\text{snap}} \|_a \leq \left\| (I-\tilde{\pi}_0)u-u_{\text{snap}} \right\|_a + \left\| \tilde{\pi}_0 u \right\|_a \leq 2\left\| \tilde{\pi}_0 u \right\|_a.
\end{equation}
It remains to estimate $\left\| \tilde{\pi}_0 u \right\|_a$. 
Since $(I-\tilde{\pi}_0)u \in V_{\text{snap}}$, we can uniquely write it into a linear combination 
\[
(I-\tilde{\pi}_0)u = 
\sum_{E_\ell \in \mathcal{E}^H} \sum_{j=1}^{J_\ell} 
c_{j}^{(\ell)} \psi_{\text{snap}}^{\ell, j}.
\]
Then we define $\tilde{p}_{\text{snap}} \in Q_0$ by 
\[
\tilde{p}_{\text{snap}} = 
\sum_{E_\ell \in \mathcal{E}^H} \sum_{j=1}^{J_\ell} 
c_{j}^{(\ell)} p_{\text{snap}}^{\ell, j}.
\]
As a result of \eqref{eqn:def-snap}, for any $v \in V_0$, we have 
$a\left((I-\tilde{\pi}_0)u, v\right) - b(v, \tilde{p}_{\text{snap}}) = 0$. 
Together with the first equation of \eqref{eqn:model_var}, this implies  
$a(\tilde{\pi}_0 u, v) - b(v, p - \tilde{p}_{\text{snap}}) = 0$.
Moreover, if $v \in \tilde{V}_0$, using Lemma~\ref{lemma:ortho-1}, we have 
$b(v, \pi_{Q_H}p) = 0$, which allows us to write 
\begin{equation}\label{eqn:loc-1}
a(\tilde{\pi}_0 u, v) - b(v, (I-\pi_{Q_H})p - \tilde{p}_{\text{snap}}) = 0, 
\end{equation}
On the other hand, using the results from Lemma~\ref{lemma:ortho-2}
for any $K \in \mathcal{T}^H$ and $q \in Q_0(K)$, we have 
$b((I-\tilde{\pi}_0) u, q) = 0$ and $(\pi_{Q_H} f, q)_K = 0$, which allows us to rewrite 
the second equation of \eqref{eqn:model_var} as 
\begin{equation}\label{eqn:loc-2}
b(\tilde{\pi}_0 u, q) = ((I - \pi_{Q_H}) f, q)_K. 
\end{equation}
In other words, the restrictions on $K$, i.e. 
$(\tilde{\pi}_0 u, (I-\pi_{Q_H})p - \tilde{p}_{\text{snap}}) \in \tilde{V}_0(K) \times Q_0(K)$, 
satisfies the system 
\begin{equation}
\label{eqn:loc}
\begin{aligned}
a(\tilde{\pi}_0 u,v) -b(v,(I-\pi_{Q_H})p - \tilde{p}_{\text{snap}}) & = 0 & \forall v\in \tilde{V}_{0}(K), \\
b(\tilde{\pi}_0 u,q) & = ((I - \pi_{Q_H}) f, q)_K & \forall q\in Q_{0}(K).
\end{aligned}
\end{equation}
Using Lemma~\ref{lemma:loc_stab}, we conclude that 
\[
\| \tilde{\pi}_0 u \|_{a(K)}\leq C\min_{x\in K}\{\kappa(x)\}^{-1/2}H\| (I - \pi_{Q_H})f \|_{L^{2}(K)}.
\]
The desired result follows directly from \eqref{eqn:main-stab}. 
\end{proof}
\end{lemma}

The last step is to analyze the error between the global multiscale solution 
and the localized multiscale solution. 
It is obvious that the error depends on the convergence 
of the modified Richardson iterations in the construction of 
localized multiscale basis functions. 
As we will show in Lemma~\ref{lemma:glo-loc}, 
the minimum eigenvalue $\mu_{\text{min}}$ and 
the maximum eigenvalue $\mu_{\text{max}}$ 
of the matrix $A_{\text{pre}}^{-1} A$
will put a sufficient condition 
on the scalar parameter $\tau$ for convergence, 
as well as control the convergence rate. 
The following lemma estimates the eigenvalues of the matrix $A_{\text{pre}}^{-1}A$ 
with bounds related to the mesh and the discretization. 
\begin{lemma}\label{lemma:precond-eig}
Denote by $\mu_{\text{min}}$ the minimum eigenvalue and 
$\mu_{\text{max}}$ the maximum eigenvalue of the matrix $A_{\text{pre}}^{-1} A$. 
We have 
\begin{equation}\label{eqn:precond-eig}
\Lambda M_1^{-1} \leq \mu_{\text{min}} \leq \mu_{\text{max}} \leq M_2,
\end{equation}
where 
\[ \begin{split}
M_1 & = \max_{K \in \mathcal{T}^H} \left\vert \left\{ E_\ell \in \mathcal{E}^H: K \subseteq \omega_\ell\right\} \right\vert, \\
M_2 & = \max_{\ell \in \{1, \cdots, N_e\}} \left\vert \left\{ s \in \{1, \cdots, N_e\}: \tilde{V}_{\text{sms}}^{(s)} \not\perp_a \tilde{V}_{\text{sms}}^{(\ell)} \right\}\right\vert, \\
\Lambda &= \min_{\ell \in \{1, \cdots, N_e\}} \lambda_{\mathcal{J}_{\ell,2} +1}^{(\ell)}, 
\end{split} \]
and $\{\lambda_j^{(\ell)}\}$ are the eigenvalues obtained from \eqref{eqn:spectral-h}. 
\begin{proof}
Let $(\mu, \tilde{\Psi})$ be an eigenpair of the matrix $A_{\text{pre}}^{-1}A$. 
Then we have 
\[
\mu = \dfrac{\tilde{\Psi}^T A \tilde{\Psi}}{\tilde{\Psi}^T A_{\text{pre}} \tilde{\Psi}}.
\]
We denote the element in $\tilde V_{\text{sms}}$ with coefficients $\tilde{\Psi}$ by $\tilde \psi$, 
and write $\tilde{\psi} = \sum_{E_\ell \in \mathcal{E}^H} \tilde{\psi}^{(\ell)}$, 
where $\tilde{\psi}^{(\ell)} \in \tilde{V}_{\text{sms}}^{(\ell)}$. 
By the definition of the matrices $A$ and $A_{\text{pre}}$, we have 
\[
\begin{split}
\tilde{\Psi}^T A \tilde{\Psi} & = a(\tilde{\psi},\tilde{\psi}), \\
\tilde{\Psi}^T A_{\text{pre}} \tilde{\Psi} & = \sum_{E_\ell \in \mathcal{E}^H} a(\tilde{\psi}^{(\ell)},\tilde{\psi}^{(\ell)}), 
\end{split}
\] 
To obtain the lower bound of $\mu$, we note that for each $E_\ell \in \mathcal{E}^H$, 
since $\mathcal{H}_\ell(\tilde{\psi}^{(\ell)}) = \tilde{\psi}^{(\ell)} = \tilde{\psi}$ on $E_\ell$, we have 
\[
\begin{split}
a(\tilde{\psi}^{(\ell)},\tilde{\psi}^{(\ell)}) 
& = \int_{\omega_\ell} \kappa^{-1} \left\vert \tilde{\psi}^{(\ell)} \right\vert^2 \\
& \leq \left(\lambda_{\mathcal{J}_{\ell,2} +1}^{(\ell)}\right)^{-1} \int_{\omega_\ell} \kappa^{-1} \left\vert \mathcal{H}(\tilde{\psi}^{(\ell)}) \right\vert^2 \\ 
& \leq \left(\lambda_{\mathcal{J}_{\ell,2} +1}^{(\ell)}\right)^{-1} \int_{\omega_\ell} \kappa^{-1} \left\vert \tilde{\psi} \right\vert^2 \\
& \leq \Lambda^{-1} \sum_{K \subseteq \omega_\ell} \int_{K} \kappa^{-1} \left\vert \tilde{\psi} \right\vert^2. 
\end{split}
\]
Summing over $E_\ell \in \mathcal{E}^H$, we have 
\[
\begin{split}
\sum_{E_\ell \in \mathcal{E}^H} a(\tilde{\psi}^{(\ell)},\tilde{\psi}^{(\ell)}) 
& \leq \Lambda^{-1} \sum_{E_\ell \in \mathcal{E}^H} \sum_{K \subseteq \omega_\ell} \int_{K} \kappa^{-1} \left\vert \tilde{\psi} \right\vert^2 \\
& \leq M_1 \Lambda^{-1} \int_{D} \kappa^{-1} \left\vert \tilde{\psi} \right\vert^2 \\
& = M_1 \Lambda^{-1} a(\tilde{\psi}, \tilde{\psi}). 
\end{split}
\]
On the other hand, to obtain the upper bound of $\mu$, using Cauchy-Schwarz inequality, we have 
\[
\begin{split}
a(\tilde{\psi}, \tilde{\psi})
& = a\left(\sum_{E_\ell \in \mathcal{E}^H} \tilde{\psi}^{(\ell)}, \sum_{E_s \in \mathcal{E}^H} \tilde{\psi}^{(s)} \right) \\
& = \sum_{E_\ell \in \mathcal{E}^H}  \sum_{E_s \in \mathcal{E}^H} a(\tilde{\psi}^{(\ell)}, \tilde{\psi}^{(s)}) \\
& \leq M_2 \sum_{E_\ell \in \mathcal{E}^H} a(\tilde{\psi}^{(\ell)}, \tilde{\psi}^{(\ell)}).
\end{split} 
\]
This completes the proof. 
\end{proof}
\end{lemma}

The following lemma provides an estimate 
for the error between the global multiscale solution 
and the local solution, 
given that the scalar parameter $\tau$ in the iterative scheme 
is sufficiently small and controlled by 
the minimum eigenvalue $\mu_{\text{min}}$ and 
the maximum eigenvalue $\mu_{\text{max}}$ 
of the matrix $A_{\text{pre}}^{-1} A$. 
\begin{lemma}\label{lemma:glo-loc}
Let $(u_{\text{glo}},p_{\text{glo}})$ be the solution in \eqref{eqn:glo_msm} and 
$(u_{\text{ms}},p_{\text{ms}})$ be the solution in \eqref{eqn:loc_msm}. 
Suppose $0< \tau \leq 2(\mu_{\text{min}}+\mu_{\text{max}})^{-1}$.
We have
\begin{equation}
\| u_{\text{glo}} - u_{\text{ms}} \|_a \leq  \exp\left(- \dfrac{k \theta \mu_{\text{min}}}{\mu_{\text{max}}}\right) \| (I - \tilde{\pi}_{\text{sms}}) u_{\text{glo}} \|_a,
\end{equation}
where $k$ is the number of iterations in the construction of multiscale basis functions in \eqref{eqn:loc_basis} and $$\theta := \displaystyle{\frac{\tau (\mu_{\text{min}}+\mu_{\text{max}})}{2}} \in (0,1].$$ 
\begin{proof}
From Lemma~\ref{lemma:snap-glo}, we have $u_{\text{glo}} = u_{\text{snap}}$. 
Subtracting \eqref{eqn:loc_msm} from \eqref{eqn:snapshot_solution}, we deduce 
\begin{equation}
\label{eqn:snap-loc}
\begin{aligned}
a(u_{\text{snap}}-u_{\text{ms}},v) -b(v,p_{\text{snap}}-p_{\text{ms}}) & = 0 & \forall v\in V_{\text{ms}}, \\
b(u_{\text{snap}}-u_{\text{ms}},q) & =0 & \forall q\in Q_{H}.
\end{aligned}
\end{equation}
We note that $u_{\text{snap}}-u_{\text{ms}} \in V_{\text{snap}}$. 
Using that fact that $\nabla \cdot V_{\text{snap}} \subset Q_H$, we take 
$q = \nabla \cdot (u_{\text{snap}}-u_{\text{ms}}) \in Q_H$ 
in the second equation of \eqref{eqn:snap-loc} to deduce that 
$\nabla \cdot \left(u_{\text{snap}}-u_{\text{ms}}\right)= 0$ almost everywhere in $D$, 
and therefore 
\[
b(u_{\text{snap}}-u_{\text{ms}}, p - p_{\text{snap}}) 
= \int_{D} (p - p_{\text{snap}}) \nabla \cdot \left(u_{\text{snap}}-u_{\text{ms}}\right) = 0.
\]
Since $u_{\text{snap}} = u_{\text{glo}} \in V_{\text{glo}}$, one can uniquely write it into a linear combination such that 
\[
u_{\text{snap}} = 
\sum_{E_\ell \in \mathcal{E}^H} \sum_{j=1}^{\mathcal{J}_\ell} 
c_{j}^{(\ell)} \psi_{\text{ms},j}^{(\ell)}.
\]
We denote $\hat{u}_{\text{snap}} \in V_{\text{ms}}$ by 
\[
\hat{u}_{\text{snap}} = 
\sum_{E_\ell \in \mathcal{E}^H} \sum_{j=1}^{\mathcal{J}_\ell} 
c_{j}^{(\ell)} \psi_{\text{ms},j}^{(\ell),k}.
\]
Then we have 
\[
\begin{split}
u_{\text{snap}} - \hat{u}_{\text{snap}} = \sum_{E_\ell \in \mathcal{E}^H} \sum_{j=1}^{\mathcal{J}_\ell} 
c_{j}^{(\ell)} (\tilde{\psi}_{\text{ms},j}^{(\ell)} - \tilde{\psi}_{\text{ms},j}^{(\ell),k}) \in \tilde{V}_{\text{sms}} \subset V_{\text{snap},2}.
\end{split}
\]
Using Lemma~\ref{lemma:ortho-1}, we have 
$b(u_{\text{snap}} - \hat{u}_{\text{snap}}, p - p_{\text{snap}}) = 0$, and therefore 
\[ b(u_{\text{ms}} - \hat{u}_{\text{snap}}, p - p_{\text{snap}}) = 0.
\] 
Taking $v = u_{\text{ms}} - \hat{u}_{\text{snap}} \in V_{\text{ms}}$ in the first equation of \eqref{eqn:snap-loc}, we have 
\[ a(u_{\text{snap}}-u_{\text{ms}},u_{\text{snap}}-u_{\text{ms}}) = a(u_{\text{snap}}-u_{\text{ms}}, u_{\text{snap}} - \hat{u}_{\text{snap}}), \]
which implies $\| u_{\text{snap}}-u_{\text{ms}} \|_a \leq \| u_{\text{snap}} - \hat{u}_{\text{snap}} \|_a$. 
In order to estimate $\| u_{\text{snap}} - \hat{u}_{\text{snap}} \|_a$, we introduce the Cholesky factorization 
$A = LL^T$ of the symmetric and positive definite matrix $A$, and define 
\[ e^k := \sum_{E_\ell \in \mathcal{E}^H} \sum_{j=1}^{\mathcal{J}_\ell} c_j^{(\ell)} L^T(\tilde\Psi_{\text{ms}, j}^{(\ell),k} - \tilde\Psi_{\text{ms}, j}^{(\ell)}) \in \mathbb{R}^\mathcal{M}. \]
Then it is straightforward to see that 
\[ 
\|e^k\|_2 = \| u_{\text{snap}} - \hat{u}_{\text{snap}} \|_a,
\]
where $\| \tilde{\Psi} \|_2$ denotes the Euclidean norm of a vector $\tilde{\Psi} \in \mathbb{R}^{\mathcal{M}}$. 
Combining \eqref{eqn:main-iterative-mat} and \eqref{eqn:iterative-sch-3}, we observe that 
\[
\tilde\Psi_{\text{ms}, j}^{(\ell),k} - \tilde\Psi_{\text{ms}, j}^{(\ell)} = (\mathbf{I}_{\mathcal{M}} - \tau A_{\text{pre}}^{-1} A) (\tilde\Psi_{\text{ms}, j}^{(\ell), k-1} - \tilde\Psi_{\text{ms}, j}^{(\ell)}). 
\]
Multiplying $L^T$ on the left and taking the linear combination with coefficients $c_j^{(\ell)}$, we have 
\[
e^k = (\mathbf{I}_{\mathcal{M}} - \tau L^T A_{\text{pre}}^{-1} L) e^{k-1}. 
\]
By induction, we have 
\[
e^k = (\mathbf{I}_{\mathcal{M}} - \tau L^T A_{\text{pre}}^{-1} L)^k e^0. 
\]
Taking Euclidean norm on both sides, we have 
\[
\| u_{\text{snap}} - \hat{u}_{\text{snap}} \|_a \leq \left ( \| \mathbf{I}_{\mathcal{M}} - \tau L^T A_{\text{pre}}^{-1} L \|_2\right)^k \| \tilde{\pi}_{\text{sms}} u_{\text{snap}} \|_a. 
\]
Since $\tilde{\pi}_{\text{sms}} u_{\text{snap}} \in \tilde{V}_{\text{sms}} \subset V_{\text{snap},2}$, 
by Lemma~\ref{lemma:ortho-1}, we have $b(\tilde{\pi}_{\text{sms}} u_{\text{snap}}, p_{\text{glo}}) = 0$. 
Taking $v = \tilde{\pi}_{\text{sms}} u_{\text{snap}} \in V_{\text{snap}}$ in the first equation of \eqref{eqn:snapshot_solution}, we have 
\[ a(\tilde{\pi}_{\text{sms}} u_{\text{snap}}, \tilde{\pi}_{\text{sms}} u_{\text{snap}}) 
= -a((I-\tilde{\pi}_{\text{sms}}) u_{\text{snap}}, \tilde{\pi}_{\text{sms}} u_{\text{snap}}), \]
which implies $\| \tilde{\pi}_{\text{sms}} u_{\text{snap}} \|_a \leq \| (I - \tilde{\pi}_{\text{sms}}) u_{\text{snap}} \|_a$. 
It remains to estimate the spectral norm $\| \mathbf{I}_{\mathcal{M}} - \tau L^T A_{\text{pre}}^{-1} L \|_2$. 
We note that $(\nu, \tilde{\Psi})$ is an eigenpair of the matrix $\mathbf{I}_{\mathcal{M}} - \tau L^T A_{\text{pre}}^{-1} L$ 
if and only if $(\tau^{-1}(1- \nu), L^{-T} \tilde{\Psi})$ is an eigenpair of the matrix $A_{\text{pre}}^{-1} A$. 
Therefore, by the assumption on the scalar parameter $\tau$, we have  $\tau < 2\mu_{\text{max}}^{-1}$ and 
\[ -1 < 1 - \mu_{\text{max}} \tau \leq \nu \leq 1 - \mu_{\text{min}} \tau < 1, \] 
which implies $\| \mathbf{I}_{\mathcal{M}} - \tau L^T A_{\text{pre}}^{-1} L \|_2 < 1$. 
Moreover, if we take $\tau \leq 2(\mu_{\text{min}}+\mu_{\text{max}})^{-1}$, then we have 
\[ 
\| \mathbf{I}_{\mathcal{M}} - \tau L^T A_{\text{pre}}^{-1} L \|_2 
= 1 - \tau \mu_{\text{min}} 
\leq 1 - \dfrac{\theta\mu_{\text{min}}}{\mu_{\text{max}}} 
\leq \exp\left(- \dfrac{\theta \mu_{\text{min}}}{\mu_{\text{max}}} \right). 
\]
This completes the proof. 
\end{proof}
\end{lemma}

Finally, we present a sufficient condition for linear convergence. 
\begin{theorem} \label{thm:main}
Let $(u,p)$ be the solution in \eqref{eqn:model_var} and 
$(u_{\text{ms}},p_{\text{ms}})$ be the solution in \eqref{eqn:loc_msm}. 
Suppose $\tau \leq M_2^{-1}$ 
and $k \geq \tau^{-1} \Lambda^{-2} M_1^2 M_2 \log(H^{-1})$.  
We have
\begin{equation}
\| u - u_{\text{ms}} \|_a \leq H\left(
C\min_{x\in\Omega}\{\kappa(x)\}^{-1/2}\|(I - \pi_{Q_H}) f\|_{L^{2}} + \| (I - \tilde{\pi}_{\text{sms}}) u_{\text{glo}} \|_a\right).
\end{equation}
\begin{proof}
By Lemma~\ref{lemma:precond-eig}, 
$M_1^{-1} \Lambda \leq \mu_{\text{min}} \leq \mu_{\text{max}} \leq M_2$. 
Since $\tau \leq M_2^{-1} \leq 2(\mu_{\text{min}}+\mu_{\text{max}})^{-1}$, by Lemma~\ref{lemma:glo-loc}, we have 
\[
\| u_{\text{glo}} - u_{\text{ms}} \|_a \leq  \exp\left(- \dfrac{k \theta \mu_{\text{min}}}{\mu_{\text{max}}}\right) \| (I - \tilde{\pi}_{\text{sms}}) u_{\text{glo}} \|_a.
\]
Moreover, we have 
\[ 
\dfrac{\theta \mu_{\text{min}}}{\mu_{\text{max}}} 
= \dfrac{\tau \mu_{\text{min}}(\mu_{\text{min}} + \mu_{\text{max}})}{2\mu_{\text{max}}} 
\geq \dfrac{\tau \mu_{\text{min}}^2}{\mu_{\text{max}}} 
\geq \dfrac{\tau \Lambda^2}{M_1^2 M_2}.
\]
With the assumption  $k \geq \tau^{-1} \Lambda^{-2} M_1^2 M_2 \log(H^{-1})$, we have 
\[ 
\exp\left(- \dfrac{k \theta \mu_{\text{min}}}{\mu_{\text{max}}}\right) 
\leq \exp\left(- \dfrac{k \tau \Lambda^2}{M_1^2 M_2}\right) \leq H. 
\]
Using a triangle inequality and invoking the result from Lemma~\ref{lemma:weak-glo}, we obtain the desired result. 
\end{proof}
\end{theorem}

\section{Numerical experiments}\label{sec:numerics}
In this section, we provide some numerical results to demonstrate the efficiency of the proposed iterative multiscale construction. 
We set the computational domain to be $D = (0,1)^2$. We use a rectangular 
mesh for the partition of the domain dividing $D$ into several coarse square elements to obtain a coarse grid $\mathcal{T}^H$ with mesh size $H>0$. Further, we divide each coarse element into several fine square elements such that the overall fine resolution is 
$256 \times 256$ with fine mesh size $h = \sqrt{2}/256$. We refer this partition to be a fine grid $\mathcal{T}^h$. 
The reference solution $(u,p)$ is solved on this fine grid by the lowest order Raviart-Thomas element ($RT0$).
In the following, we define $L^2$ error of the pressure variable and  the energy error of velocity variable as follows:
$$ e_{2} := \frac{\norm{p - p_{\text{ms}}}_{L^2(D)}}{\norm{p}_{L^2(D)}} \quad \text{and} \quad 
e_a := \frac{\norm{u - u_{\text{ms}}}_{a}}{\norm{u}_{a}}.$$
Here, $(u_{\text{ms}}, p_{\text{ms}})$ is the multiscale solution obtained by solving \eqref{eqn:loc_msm} using the iteration-constructed multiscale basis functions. 
It is remarkable that under this setting of coarse mesh and the $RT0$ element that is used, we have $M_1 = 4$ and $M_2 = 3$, where $M_1$ and $M_2$ are defined in Lemma \ref{lemma:precond-eig}. 

In all the examples below, we set the initial condition $\tilde \psi_{\text{ms}, j}^{(\ell),0} = 0$. The regularization parameter is either set to be $\tau = M_2^{-1} = 1/3$ or $\tau = \tau_{\text{opt}} := 2(\mu_{\min} + \mu_{\max})^{-1}$, where $ \mu_{\min}$ and $\mu_{\max}$ are the smallest and largest eigenvalues of the matrix $A_{\text{pre}}^{-1} A$, respectively. 
We denote $k \in \mathbb{N}$ the number of iteration level. 
We remark that the choice of the regularization parameter $\tau$ is crucial in the proposed iterative scheme for multiscale basis functions. In practice, one may not have any a priori information about the spectrum of the matrix $A_{\text{pre}}^{-1} A$, and the case of $\tau = M_2^{-1}$ serves as a baseline of the performance using the iterative construction. 

\begin{exmp} \label{exp1}
In this example, we consider the heterogeneous media $\kappa$ to be defined as follows: 
$$\kappa (x_1, x_2) =\cfrac{2+\sin(11\pi x_1)\sin(13\pi x_2)}{1.4+\cos(12\pi x_1)\cos(7\pi x_2)}$$
for any $(x_1, x_2) \in D$. The source function in this example is defined to be 
$$ f(x_1, x_2) = \left \{ \begin{array}{cl}
1 & \text{for} ~ x_1 \in [0, 1/2),\\
-1 & \text{for} ~x_1 \in [1/2, 1].
\end{array} \right . $$
We remark that the source function satisfies the compatibility condition $\int_D f ~dx= 0$. 
We choose $\mathcal{J}_{\ell} = 2$ to form the local multiscale space $V_{\text{snap},1}^{(\ell)} \oplus V_{\text{sms}}^{(\ell)}$ for each coarse neighborhood $\omega_{\ell}$. 
In Tables \ref{tab:exp1-v-subopt} and \ref{tab:exp1-p-subopt}, we show the velocity and pressure errors with $\tau = M_2^{-1}$.
The results of errors using $\tau = \tau_{\text{opt}}$ are reported in Tables \ref{tab:exp1-v} and \ref{tab:exp1-p}. 
In both the cases, one can observe the convergence with respect to the coarse mesh size $H$ and the number of iteration $k$ for constructing the multiscale basis functions while the case with $\tau = \tau_{\text{opt}}$ gives a rapider decay of energy error. It can be observed that for a given fixed coarse mesh, the iterative process for multiscale basis functions barely improve the accuracy of approximation for pressure variable since the basis functions for pressure are identical during the iteration. 

\begin{table}[ht]
\centering
\begin{tabular}{cc||ccccccc}
\hline
\multicolumn{2}{c||}{\multirow{2}{*}{$e_a$}}   & \multicolumn{7}{c}{Iteration level $k$}           \\ \cline{3-9} 
\multicolumn{2}{c||}{} &   $0$ & $1$   & $2$   & $3$   & $4$   &$5$ & $6$ \\ 
\hline \hline
\multicolumn{1}{c|}{\multirow{3}{*}{$H$}} & 
$\sqrt{2}/8$ & $16.7132\%$ & $10.9770\%$ & $7.3665\%$ & $5.0761\%$ & $3.5936\%$ & $2.6093\%$ & $1.9377\%$ \\

\multicolumn{1}{c|}{} &  
$\sqrt{2}/16$ & $7.0126\%$ & $4.7449\%$ & $3.2317\%$ & $2.2197\%$ & $1.5399\%$ & $1.0806\%$ & $0.7680\%$ \\

\multicolumn{1}{c|}{} &  
$\sqrt{2}/32$ & $2.4646\%$ & $1.6716\%$ & $1.1307\%$ & $0.7645\%$ & $0.5173\%$ & $0.3505\%$ & $0.2379\%$ \\
\hline
\end{tabular}
\caption{Energy errors with $\mathcal{J}_{\ell} = 2$, $\tau = M_2^{-1}$, and varying $k$ and $H$ (Example \ref{exp1}).}
\label{tab:exp1-v-subopt}
\end{table}

\begin{table}[ht]
\centering
\begin{tabular}{cc||ccccccc}
\hline
\multicolumn{2}{c||}{\multirow{2}{*}{$e_2$}}   & \multicolumn{7}{c}{Iteration level $k$}           \\ \cline{3-9} 
\multicolumn{2}{c||}{} &   $0$ & $1$   & $2$   & $3$   & $4$   &$5$ & $6$ \\ 
\hline \hline
\multicolumn{1}{c|}{\multirow{3}{*}{$H$}} & 
$\sqrt{2}/8$ & $25.8797\%$ & $25.7675\%$ & $25.7472\%$ & $25.7432\%$ & $25.7423\%$ & $25.7421\%$ & $25.7421\%$ \\

\multicolumn{1}{c|}{} &  
$\sqrt{2}/16$ & $13.3372\%$ & $13.3304\%$ & $13.3289\%$ & $13.3286\%$ & $13.3286\%$ & $13.3286\%$ & $13.3285\%$ \\

\multicolumn{1}{c|}{} &  
$\sqrt{2}/32$ & $6.7420\%$ & $6.7418\%$ & $6.7418\%$ & $6.7417\%$ & $6.7417\%$ & $6.7417\%$ & $6.7417\%$ \\

\hline
\end{tabular}
\caption{$L^2$ errors with $\mathcal{J}_{\ell} = 2$, $\tau = M_2^{-1}$, and varying $k$ and $H$ (Example \ref{exp1}).}
\label{tab:exp1-p-subopt}
\end{table}

\begin{table}[ht]
\centering
\begin{tabular}{cc||ccccc}
\hline
\multicolumn{2}{c||}{\multirow{2}{*}{$e_a$}}   & \multicolumn{5}{c}{Iteration level $k$}           \\ \cline{3-7} 
\multicolumn{2}{c||}{} &   $0$ & $1$   & $2$   & $3$   & $4$   \\ 
\hline \hline
\multicolumn{1}{c|}{\multirow{3}{*}{$H$}} & 
$\sqrt{2}/8$ & $16.7132\%$ & $3.8980\%$ & $1.4137\%$ & $0.5778\%$ & $0.2575\%$ \\

\multicolumn{1}{c|}{} &  
$\sqrt{2}/16$ & $7.0126\%$ & $1.1017\%$ & $0.3679\%$ & $0.1401\%$ & $0.0563\%$ \\

\multicolumn{1}{c|}{} &  
$\sqrt{2}/32$ & $2.4646\%$ & $0.1694\%$ & $0.0253\%$ & $0.0052\%$ & $0.0013\%$ \\
\hline
\end{tabular}
\caption{Energy errors with $\mathcal{J}_{\ell} = 2$, $\tau = \tau_{\text{opt}}$, and varying $k$ and $H$ (Example \ref{exp1}).}
\label{tab:exp1-v}
\end{table}

\begin{table}[ht]
\centering
\begin{tabular}{cc||ccccc}
\hline
\multicolumn{2}{c||}{\multirow{2}{*}{$e_2$}}   & \multicolumn{5}{c}{Iteration level $k$}           \\ \cline{3-7} 
\multicolumn{2}{c||}{} &   $0$ & $1$   & $2$   & $3$   & $4$   \\ 
\hline \hline
\multicolumn{1}{c|}{\multirow{3}{*}{$H$}} & 
$\sqrt{2}/8$ & $25.8797\%$ & $25.7425\%$ & $25.7420\%$ & $25.7420\%$ & $25.7420\%$ \\

\multicolumn{1}{c|}{} &  
$\sqrt{2}/16$ & $13.3372\%$ & $13.3286\%$ & $13.3285\%$ & $13.3285\%$ & $13.3285\%$ \\

\multicolumn{1}{c|}{} &  
$\sqrt{2}/32$ & $6.7420\%$ & $6.7417\%$ & $6.7417\%$ & $6.7417\%$ & $6.7417\%$ \\

\hline
\end{tabular}
\caption{$L^2$ errors with $\mathcal{J}_{\ell} = 2$, $\tau = \tau_{\text{opt}}$, and varying $k$ and $H$ (Example \ref{exp1}).}
\label{tab:exp1-p}
\end{table}
\end{exmp}

\begin{exmp} \label{exp2}
In this example, we consider a permeability field which is of high value of contrast. The permeability field is depicted in Figure \ref{fig:kappa} (left). The source function is defined to be a piecewise constant function satisfying 
$$ f(x_1, x_2) = \left \{ \begin{array}{cl}
1 & \text{for} ~ (x_1, x_2) \in [7/8, 1] \times [1/2, 5/8], \\
-1 & \text{for} ~ (x_1, x_2) \in [7/8, 1] \times [7/8, 1], \\
0 & \text{otherwise}.
\end{array} \right . $$
The source function satisfies the compatibility condition. The errors in the case with $\tau = M_2^{-1}$ are shown in Tables \ref{tab:exp2-v-subopt} and \ref{tab:exp2-p-subopt} and the results with $\tau = \tau_{\text{opt}}$ are recorded in Tables \ref{tab:exp2-v} and \ref{tab:exp2-p}. 

The iterative construction help enhance the accuracy of the approximation of the velocity variable. Comparing to the non-iterative CEM basis construction in \cite{chung2018constraintmixed}, the iterative construction of the basis functions provides a flexible approach to compute the basis functions within a desire threshold of accuracy. Moreover, the iterative approach of constructing basis functions has smaller marginal computational cost from $k$-th level iteration to $k+1$-th level iteration given a fixed coarse grid; while decreasing the coarse mesh size requires more computation. 

\begin{figure}[ht!]
\centering
\includegraphics[width=2.7in]{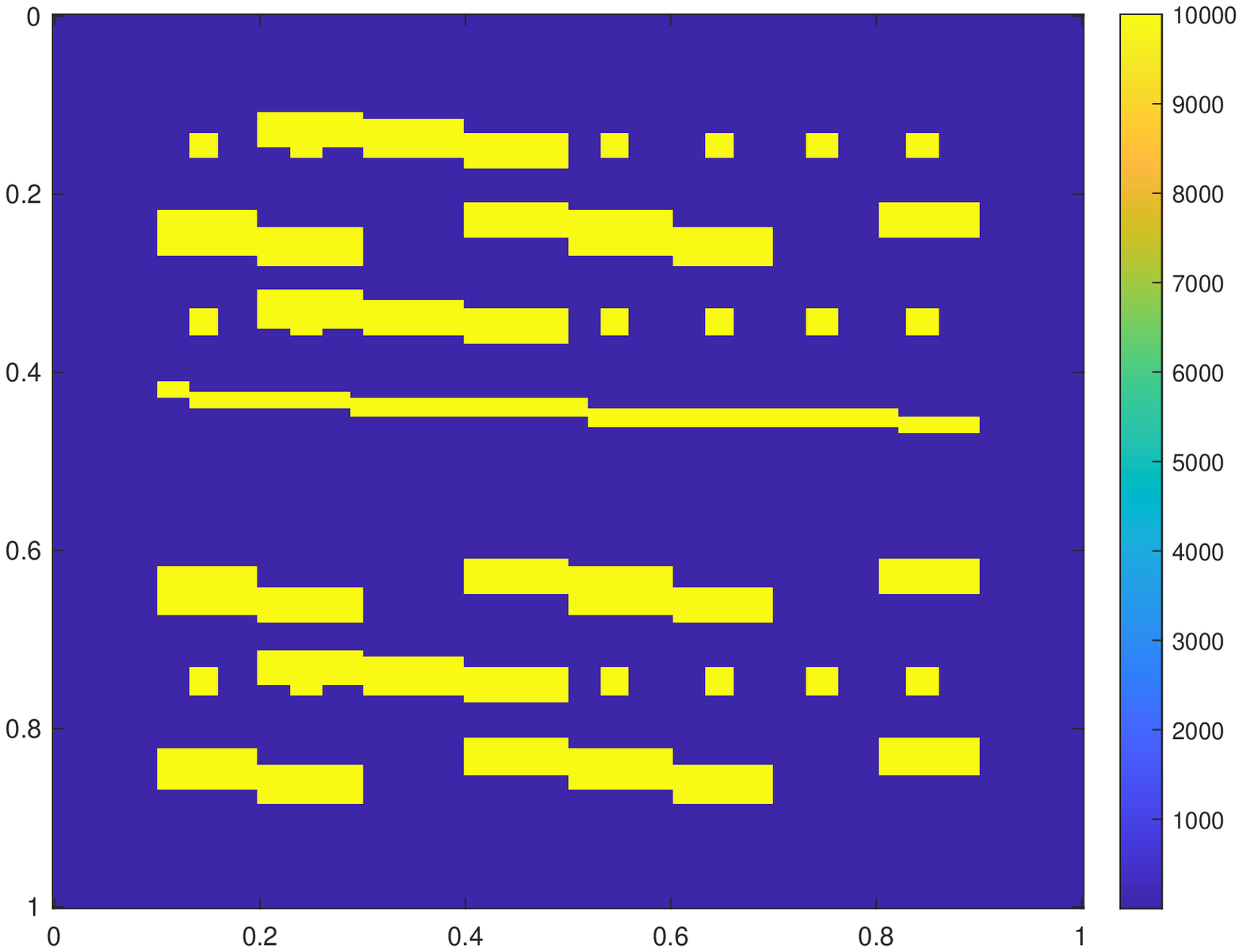}
\qquad
\includegraphics[width=2.7in]{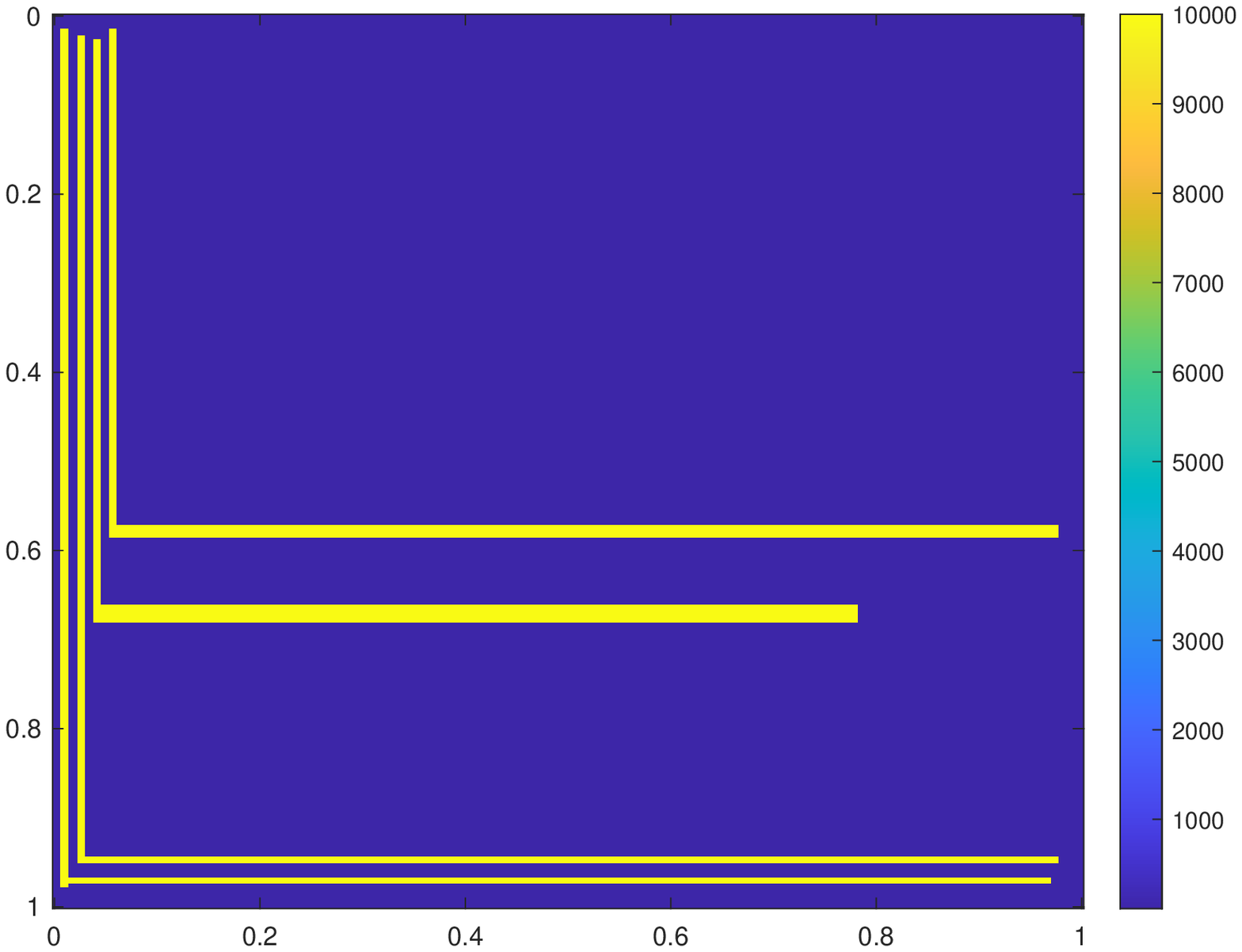}
\caption{Permeability fields. Left: Example \ref{exp2}. Right: Example \ref{exp3}.} 
\label{fig:kappa}
\end{figure}

\begin{table}[ht]
\centering
\begin{tabular}{cc||ccccccc}
\hline
\multicolumn{2}{c||}{\multirow{2}{*}{$e_a$}}   & \multicolumn{7}{c}{Iteration level $k$}           \\ \cline{3-9} 
\multicolumn{2}{c||}{} &   $0$ & $1$   & $2$   & $3$   & $4$   &$5$ & $6$ \\ 
\hline \hline
\multicolumn{1}{c|}{\multirow{3}{*}{$H$}} & 
 $\sqrt{2}/8$ & $28.5514\%$ & $19.8872\%$ & $13.9714\%$ & $9.7932\%$ & $6.8657\%$ & $4.8412\%$ & $3.4483\%$ \\
 
\multicolumn{1}{c|}{} &  
$\sqrt{2}/16$ & $21.4083\%$ & $16.2509\%$ & $11.9394\%$ & $8.5096\%$ & $5.9568\%$ & $4.1415\%$ & $2.8803\%$ \\

\multicolumn{1}{c|}{} &  
$\sqrt{2}/32$ & $16.8212\%$ & $12.8958\%$ & $9.3932\%$ & $6.5901\%$ & $4.5287\%$ & $3.0844\%$ & $2.0960\%$ \\

\hline
\end{tabular}
\caption{Energy errors with $\mathcal{J}_{\ell} = 2$, $\tau = M_2^{-1}$, and varying $k$ and $H$ (Example \ref{exp2}).}
\label{tab:exp2-v-subopt}
\end{table}

\begin{table}[ht]
\centering
\begin{tabular}{cc||ccccccc}
\hline
\multicolumn{2}{c||}{\multirow{2}{*}{$e_2$}}   & \multicolumn{7}{c}{Iteration level $k$}           \\ \cline{3-9} 
\multicolumn{2}{c||}{} &   $0$ & $1$   & $2$   & $3$   & $4$   &$5$ & $6$ \\ 
\hline \hline
\multicolumn{1}{c|}{\multirow{3}{*}{$H$}} & 
$\sqrt{2}/8$ & $40.4620\%$ & $38.3771\%$ & $37.3959\%$ & $37.0434\%$ & $36.9396\%$ & $36.9119\%$ & $36.9048\%$ \\

\multicolumn{1}{c|}{} &  
$\sqrt{2}/16$ & $30.6075\%$ & $27.8506\%$ & $26.3574\%$ & $25.7909\%$ & $25.6255\%$ & $25.5838\%$ & $25.5738\%$ \\

\multicolumn{1}{c|}{} &  
$\sqrt{2}/32$ & $19.8788\%$ & $16.3759\%$ & $14.5293\%$ & $13.8820\%$ & $13.7100\%$ & $13.6702\%$ & $13.6616\%$ \\

\hline
\end{tabular}
\caption{$L^2$ errors with $\mathcal{J}_{\ell} = 2$, $\tau = M_2^{-1}$, and varying $k$ and $H$ (Example \ref{exp2}).}
\label{tab:exp2-p-subopt}
\end{table}

\begin{table}[ht!]
\centering
\begin{tabular}{cc||ccccc}
\hline
\multicolumn{2}{c||}{\multirow{2}{*}{$e_a$}}   & \multicolumn{5}{c}{Iteration level $k$}           \\ \cline{3-7} 
\multicolumn{2}{c||}{} &   $0$ & $1$   & $2$   & $3$   & $4$   \\ 
\hline \hline
\multicolumn{1}{c|}{\multirow{3}{*}{$H$}} & $\sqrt{2}/8$ & $28.5514\%$ & $7.6068\%$ & $3.5353\%$ & $1.8574\%$ & $1.0224\%$ \\
\multicolumn{1}{c|}{} &  $\sqrt{2}/16$ & $21.4083\%$ & $3.6325\%$ & $0.9163\%$ & $0.2764\%$ & $0.0975\%$ \\
\multicolumn{1}{c|}{} &  $\sqrt{2}/32$ & $16.8212\%$ & $1.8208\%$ & $0.2752\%$ & $0.0527\%$ & $0.0126\%$ \\
\hline
\end{tabular}
\caption{Energy errors with $\mathcal{J}_{\ell} = 2$, $\tau = \tau_{\text{opt}}$, and varying $k$ and $H$ (Example 
\ref{exp2}).}
\label{tab:exp2-v}
\end{table}

\begin{table}[ht!]
\centering
\begin{tabular}{cc||ccccc}
\hline
\multicolumn{2}{c||}{\multirow{2}{*}{$e_2$}}   & \multicolumn{5}{c}{Iteration level $k$}           \\ \cline{3-7} 
\multicolumn{2}{c||}{} &   $0$ & $1$   & $2$   & $3$   & $4$   \\ 
\hline \hline
\multicolumn{1}{c|}{\multirow{3}{*}{$H$}} & 
$\sqrt{2}/8$ & $40.4620\%$ & $36.9186\%$ & $36.9026\%$ & $36.9023\%$ & $36.9022\%$ \\
\multicolumn{1}{c|}{} & 
$\sqrt{2}/16$ & $30.6075\%$ & $25.5770\%$ & $25.5708\%$ & $25.5708\%$ & $25.5708\%$ \\
\multicolumn{1}{c|}{} &  
$\sqrt{2}/32$ & $19.8788\%$ & $13.6605\%$ & $13.6593\%$ & $13.6593\%$ & $13.6593\%$ \\
\hline
\end{tabular}
\caption{$L^2$ errors with $\mathcal{J}_{\ell} = 2$, $\tau = \tau_{\text{opt}}$, and varying $k$ and $H$ (Example \ref{exp2}).}
\label{tab:exp2-p}
\end{table}
\end{exmp}

\begin{exmp} \label{exp3}
In this example, we consider a more challenging channelized permeability field and it is sketched in Figure \ref{fig:kappa} (right). The source function is defined to be 
$$ f(x_1, x_2) = \left \{ \begin{array}{cl}
1 & \text{for} ~ (x_1, x_2) \in [0, 1/8] \times [0, 1/8], \\
-1 & \text{for} ~ (x_1, x_2) \in [7/8, 1] \times [7/8, 1], \\
0 & \text{otherwise}.
\end{array} \right . $$
The source function also satisfies the compatibility condition. We set $\mathcal{J}_{\ell} = 3$ or $\mathcal{J}_{\ell} = 4$ in this example. 
The numerical results with $\mathcal{J}_{\ell} = 3$ and $\tau = M_2^{-1}$ are depicted in Tables \ref{tab:exp3-v-subopt-2} and \ref{tab:exp3-p-subopt-2}. 
The corresponding results of errors with $\mathcal{J}_{\ell} = 3$ and $\tau = \tau_{\text{opt}}$ are presented in Tables \ref{tab:exp3-v-2} and \ref{tab:exp3-p-2}. 
The errors with $\mathcal{J}_{\ell} = 4$ and $\tau = M_2^{-1}$ are shown in Tables \ref{tab:exp3-v-subopt-3} and \ref{tab:exp3-p-subopt-3} while 
those with $\mathcal{J}_{\ell} = 4$ and $\tau = \tau_{\text{opt}}$ are shown in Tables \ref{tab:exp3-v-3} and \ref{tab:exp3-p-3}. 
We remark that in this example, one has to include more basis functions to form the auxiliary space $V_{\text{sms}}$ in order to obtain sharp convergence rate with respect to the number of iterations. 

In this example, besides the observation of decay of velocity errors, we can observe the decay of the pressure error during the iterations and the error stalls eventually at a smaller magnitude. For instance, when $H = \sqrt{2}/32$ and $\mathcal{J}_{\ell} = 3$, the $L^2$ error at the beginning is about $87.9352\%$; it is around the level of $e_2 = 23.5540\%$ after a few iterations with $\tau = M_2^{-1}$. 

\begin{table}[ht]
\centering
\begin{tabular}{cc||ccccccc}
\hline
\multicolumn{2}{c||}{\multirow{2}{*}{$e_a$}}   & \multicolumn{7}{c}{Iteration level $k$}           \\ \cline{3-9} 
\multicolumn{2}{c||}{} &   $0$ & $1$   & $2$   & $3$   & $4$   &$5$ & $6$ \\ 
\hline \hline
\multicolumn{1}{c|}{\multirow{3}{*}{$H$}} & 
$\sqrt{2}/8$ & $109.5043\%$ & $97.3436\%$ & $85.5803\%$ & $74.0380\%$ & $62.4770\%$ & $50.8968\%$ & $39.8888\%$ \\
\multicolumn{1}{c|}{} &  
 $\sqrt{2}/16$ & $103.6153\%$ & $94.0967\%$ & $83.2797\%$ & $71.3002\%$ & $58.4448\%$ & $45.5061\%$ & $33.8691\%$ \\
\multicolumn{1}{c|}{} &  
 $\sqrt{2}/32$ & $87.9352\%$ & $77.6234\%$ & $66.7290\%$ & $55.2908\%$ & $43.6862\%$ & $32.7943\%$ & $23.5540\%$ \\

\hline
\end{tabular}
\caption{Energy errors with $\mathcal{J}_{\ell} = 3$, $\tau = M_2^{-1}$, and varying $k$ and $H$ (Example \ref{exp3}).}
\label{tab:exp3-v-subopt-2}
\end{table}

\begin{table}[ht]
\centering
\begin{tabular}{cc||ccccccc}
\hline
\multicolumn{2}{c||}{\multirow{2}{*}{$e_2$}}   & \multicolumn{7}{c}{Iteration level $k$}           \\ \cline{3-9} 
\multicolumn{2}{c||}{} &   $0$ & $1$   & $2$   & $3$   & $4$   &$5$ & $6$ \\ 
\hline \hline
\multicolumn{1}{c|}{\multirow{3}{*}{$H$}} & 
 $\sqrt{2}/8$ & $100.3050\%$ & $90.5131\%$ & $81.1143\%$ & $72.4292\%$ & $64.6190\%$ & $58.3690\%$ & $54.4724\%$ \\ 

\multicolumn{1}{c|}{} &  
 $\sqrt{2}/16$ & $88.2374\%$ & $79.1280\%$ & $68.0928\%$ & $55.9858\%$ & $44.3716\%$ & $35.5116\%$ & $30.7835\%$ \\

\multicolumn{1}{c|}{} &  
 $\sqrt{2}/32$ & $80.1057\%$ & $68.4491\%$ & $54.1957\%$ & $39.5810\%$ & $28.0122\%$ & $21.0571\%$ & $17.9814\%$ \\

\hline
\end{tabular}
\caption{$L^2$ errors with $\mathcal{J}_{\ell} = 3$, $\tau = M_2^{-1}$, and varying $k$ and $H$ (Example \ref{exp3}).}
\label{tab:exp3-p-subopt-2}
\end{table}

\begin{table}[ht]
\centering
\begin{tabular}{cc||ccccc}
\hline
\multicolumn{2}{c||}{\multirow{2}{*}{$e_a$}}   & \multicolumn{5}{c}{Iteration level $k$}           \\ \cline{3-7} 
\multicolumn{2}{c||}{} &   $0$ & $1$   & $2$   & $3$   & $4$   \\ 
\hline \hline
\multicolumn{1}{c|}{\multirow{3}{*}{$H$}} & 
$\sqrt{2}/8$ & $109.5043\%$ & $47.8256\%$ & $16.3427\%$ & $6.1077\%$ & $3.0773\%$ \\

\multicolumn{1}{c|}{} &  
 $\sqrt{2}/16$ & $103.6153\%$ & $34.1147\%$ & $10.1379\%$ & $6.4740\%$ & $5.8133\%$ \\

\multicolumn{1}{c|}{} &  
 $\sqrt{2}/32$ & $87.9352\%$ & $4.2297\%$ & $0.1858\%$ & $0.0135\%$ & $0.0013\%$ \\

\hline
\end{tabular}
\caption{Energy errors with $\mathcal{J}_{\ell} = 3$, $\tau = \tau_{\text{opt}}$, and varying $k$ and $H$ (Example \ref{exp3}).}
\label{tab:exp3-v-2}
\end{table}

\begin{table}[ht]
\centering
\begin{tabular}{cc||ccccc}
\hline
\multicolumn{2}{c||}{\multirow{2}{*}{$e_2$}}   & \multicolumn{5}{c}{Iteration level $k$}           \\ \cline{3-7} 
\multicolumn{2}{c||}{} &   $0$ & $1$   & $2$   & $3$   & $4$   \\ 
\hline \hline
\multicolumn{1}{c|}{\multirow{3}{*}{$H$}} & 
$\sqrt{2}/8$ & $100.3050\%$ & $56.5662\%$ & $51.7407\%$ & $51.6529\%$ & $51.6498\%$ \\

\multicolumn{1}{c|}{} &  
 $\sqrt{2}/16$ & $88.2374\%$ & $31.6210\%$ & $28.3823\%$ & $28.3396\%$ & $28.3278\%$ \\

\multicolumn{1}{c|}{} &  
 $\sqrt{2}/32$ & $80.1057\%$ & $16.6496\%$ & $16.6481\%$ & $16.6481\%$ & $16.6481\%$ \\

\hline
\end{tabular}
\caption{$L^2$ errors with $\mathcal{J}_{\ell} = 3$, $\tau = \tau_{\text{opt}}$, and varying $k$ and $H$ (Example \ref{exp3}).}
\label{tab:exp3-p-2}
\end{table}

\begin{table}[ht]
\centering
\begin{tabular}{cc||ccccccc}
\hline
\multicolumn{2}{c||}{\multirow{2}{*}{$e_a$}}   & \multicolumn{7}{c}{Iteration level $k$}           \\ \cline{3-9} 
\multicolumn{2}{c||}{} &   $0$ & $1$   & $2$   & $3$   & $4$   &$5$ & $6$ \\ 
\hline \hline
\multicolumn{1}{c|}{\multirow{3}{*}{$H$}} & 
$\sqrt{2}/8$ & $105.8526\%$ & $94.9011\%$ & $83.4450\%$ & $71.3884\%$ & $58.4940\%$ & $45.4424\%$ & $33.6558\%$ \\
\multicolumn{1}{c|}{} &  
$\sqrt{2}/16$ & $99.2377\%$ & $90.4072\%$ & $79.8811\%$ & $67.9752\%$ & $54.9688\%$ & $41.8101\%$ & $30.1055\%$ \\
 
\multicolumn{1}{c|}{} &  
$\sqrt{2}/32$ & $70.3995\%$ & $56.9658\%$ & $43.5425\%$ & $31.6116\%$ & $22.1611\%$ & $15.2284\%$ & $10.3402\%$ \\

\hline
\end{tabular}
\caption{Energy errors with $\mathcal{J}_{\ell} = 4$, $\tau = M_2^{-1}$, and varying $k$ and $H$ (Example \ref{exp3}).}
\label{tab:exp3-v-subopt-3}
\end{table}

\begin{table}[ht]
\centering
\begin{tabular}{cc||ccccccc}
\hline
\multicolumn{2}{c||}{\multirow{2}{*}{$e_2$}}   & \multicolumn{7}{c}{Iteration level $k$}           \\ \cline{3-9} 
\multicolumn{2}{c||}{} &   $0$ & $1$   & $2$   & $3$   & $4$   &$5$ & $6$ \\ 
\hline \hline
\multicolumn{1}{c|}{\multirow{3}{*}{$H$}} & 
$\sqrt{2}/8$ & $90.4620\%$ & $83.5678\%$ & $76.7788\%$ & $69.1575\%$ & $61.3953\%$ & $55.7326\%$ & $52.9761\%$ \\
\multicolumn{1}{c|}{} &  
 $\sqrt{2}/16$ & $76.3118\%$ & $69.1737\%$ & $60.9587\%$ & $51.4130\%$ & $41.3954\%$ & $33.6981\%$ & $29.9219\%$ \\

\multicolumn{1}{c|}{} &  
 $\sqrt{2}/32$ & $53.3795\%$ & $38.6773\%$ & $26.7402\%$ & $20.0542\%$ & $17.5450\%$ & $16.8540\%$ & $16.6923\%$ \\

\hline
\end{tabular}
\caption{$L^2$ errors with $\mathcal{J}_{\ell} = 4$, $\tau = M_2^{-1}$, and varying $k$ and $H$ (Example \ref{exp3}).}
\label{tab:exp3-p-subopt-3}
\end{table}

\begin{table}[ht]
\centering
\begin{tabular}{cc||ccccc}
\hline
\multicolumn{2}{c||}{\multirow{2}{*}{$e_a$}}   & \multicolumn{5}{c}{Iteration level $k$}           \\ \cline{3-7} 
\multicolumn{2}{c||}{} &   $0$ & $1$   & $2$   & $3$   & $4$   \\ 
\hline \hline
\multicolumn{1}{c|}{\multirow{3}{*}{$H$}} & 
$\sqrt{2}/8$ & $105.8526\%$ & $25.8095\%$ & $5.2167\%$ & $1.2168\%$ & $0.2928\%$ \\
\multicolumn{1}{c|}{} &  
 $\sqrt{2}/16$ & $99.2377\%$ & $5.0349\%$ & $0.3430\%$ & $0.0314\%$ & $0.0031\%$ \\

\multicolumn{1}{c|}{} &  
 $\sqrt{2}/32$ & $70.3995\%$ & $0.2055\%$ & $0.0028\%$ & $4.989 \times 10^{-5}\%$ & $9.160 \times 10^{-7}\%$ \\

\hline
\end{tabular}
\caption{Energy errors with $\mathcal{J}_{\ell} = 4$, $\tau = \tau_{\text{opt}}$, and varying $k$ and $H$ (Example \ref{exp3}).}
\label{tab:exp3-v-3}
\end{table}

\begin{table}[ht]
\centering
\begin{tabular}{cc||ccccc}
\hline
\multicolumn{2}{c||}{\multirow{2}{*}{$e_2$}}   & \multicolumn{5}{c}{Iteration level $k$}           \\ \cline{3-7} 
\multicolumn{2}{c||}{} &   $0$ & $1$   & $2$   & $3$   & $4$   \\ 
\hline \hline
\multicolumn{1}{c|}{\multirow{3}{*}{$H$}} & 
$\sqrt{2}/8$ & $90.4620\%$ & $52.3872\%$ & $51.6513\%$ & $51.6493\%$ & $51.6493\%$ \\
\multicolumn{1}{c|}{} &  
 $\sqrt{2}/16$ & $76.3118\%$ & $28.3032\%$ & $28.3006\%$ & $28.3006\%$ & $28.3006\%$ \\

\multicolumn{1}{c|}{} &  
 $\sqrt{2}/32$ & $53.3795\%$ & $16.6481\%$ & $16.6481\%$ & $16.6481\%$ & $16.6481\%$ \\

\hline
\end{tabular}
\caption{$L^2$ errors with $\mathcal{J}_{\ell} = 4$, $\tau = \tau_{\text{opt}}$, and varying $k$ and $H$ (Example \ref{exp3}).}
\label{tab:exp3-p-3}
\end{table}
\end{exmp}

\section{Conclusion} \label{sec:conclusion}
In this work, we proposed an iterative process to construct the multiscale basis functions satisfying the property of constraint energy minimization. The procedure starts with the construction of snapshot space and we decompose the snapshot functions into the decaying and the non-decaying parts. The decaying parts are approximated iteratively via a modified Richardson scheme with an appropriate defined preconditioner, while the non-decaying parts are fixed during the iteration. 
With this set of iterative-based multiscale basis functions, we show that the energy error is of first order with respect to the coarse mesh size if sufficiently large iterations (with regularization parameter being in an appropriate range) for multiscale basis functions are conducted. Numerical experiments are provided to demonstrate the efficiency of the proposed method and confirms the theory. 

\section*{Acknowledgement}

The research of Eric Chung is partially supported by the Hong
 Kong RGC General Research Fund (Project numbers 14304719 and 14302018)
 and CUHK Faculty of Science Direct Grant 2019-20.

\bibliographystyle{plain}
\bibliography{references,references1,references2}

\end{document}